\magnification=\magstep1

\def\present#1{\noindent{\bf#1}}                     
\outer\def\give#1. {\medbreak
             \noindent{\sl#1. }}                     
\outer\def\section #1\par{\bigbreak\centerline{\S
     {\bf#1}}\nobreak\smallskip\noindent}
\def\({\left(}
\def\){\right)}
\def\contract{{\,{\vrule height5pt width0.4pt depth 0pt}
{\vrule height0.4pt width6pt depth 0pt}\,}}

\def\sqr#1#2{{\vcenter{\hrule height.#2pt              
     \hbox{\vrule width.#2pt height#1pt\kern#1pt
     \vrule width.#2pt}
     \hrule height.#2pt}}}
\def\square{\mathchoice\sqr{5.5}4\sqr{5.0}4\sqr{4.8}3\sqr{4.8}3}
\def\qed{\hskip4pt plus1fill\ $\square$\par\medbreak}



\def\cA{{\cal A}}

\def\cC{{\cal C}}
\def\cD{{\cal D}}
\def\cE{{\cal E}}
\def\cF{{\cal F}}
\def\cG{{\cal G}}
\def\cH{{\cal H}}

\def\cM{{\cal M}}
\def\cN{{\cal N}}
\def\cO{{\cal O}}

\def\cQ{{\cal Q}}
\def\cR{{\cal R}}
\def\cS{{\cal S}}

\def\cW{{\cal W}}

 
\def\CC{{\bf C}}
\def\RR{{\bf R}}

\def\C{{\bf C}}
\def\cx#1{{\CC}^{#1}}     
\def\R{{\bf R}}
\def\real#1{{\bf R}^{#1}}   

\def\bar{\overline}              
 
 
\def\bd{\partial}		 



\def\C{{\bf C}}
\def\cx#1{{\CC}^{#1}}     
 

\centerline{\bf Polynomial Diffeomorphisms of $\cx2$. IV:}
\centerline{\bf The Measure of Maximal Entropy and Laminar Currents}
 \bigskip
 \centerline{Eric Bedford\footnote{$^1$}{\rm Supported in part by NSF grant
DMS-9103585.}, Mikhail Lyubich\footnote{$^2$}{\rm Supported in part by 
NSF grant DMS-8920768 and a Sloan Research Fellowship.}, 
and John Smillie\footnote{$^3$}{\rm
Supported  in part by
 NSF grant DMS-9003101.}} \bigskip\bigskip
 
\section 1. Introduction
 
 The simplest holomorphic dynamical systems which display interesting behavior
are the polynomial maps of $\CC$. The dynamical study of
these maps began with Fatou and Julia in the 1920's and is currently a very
active area of research. If we are interested in studying {\it
invertible} holomorphic dynamical systems, then the simplest examples with
interesting behavior are probably the polynomial diffeomorphisms of $\cx2$.
These are maps $f:\cx2\to\cx2$ such that the coordinate functions of $f$ and
$f^{-1}$ are holomorphic polynomials.
 
For polynomial maps of $\CC$ the algebraic degree of the polynomial is a
useful dynamical invariant. In particular the only dynamically interesting
maps are those with degree $d$ greater than one. For polynomial
diffeomorphisms we can define the algebraic degree to be the
 maximum of the degrees of the coordinate functions. This is not, however, a
conjugacy invariant. Friedland and Milnor [FM] gave an alternative definition
of a positive integer $\deg f$ which is more natural from a dynamical point of
view. If $\deg f>1$, then $\deg f$ coincides with the minimal algebraic degree
of a diffeomorphism in the conjugacy class of $f$. As in the case of
polynomial maps of $\CC$, the polynomial diffeomorphisms $f$ with $\deg(f)=1$
are rather uninteresting. We will make the standing assumption that
$\deg(f)>1$.
 
For a polynomial map of $\cx{}$ the point at infinity is an attractor. Thus the
``recurrent'' dynamics can take place only on the set $K$ consisting of bounded
orbits.  A normal families argument shows that there is no expansion on the
interior of $K$ so ``chaotic'' dynamics can occur only on  $J=\bd K$. This set
is called the Julia set and plays a major role in the study of polynomial maps.
 
 For diffeomorphisms of $\cx2$ each of the objects $K$ and $J$ has
three analogs. Corresponding to the set $K$ in one dimension, we have the
sets $K^+$ (resp.\ $K^-$) consisting of the points whose orbits
are bounded in forward (resp.\ backward) time and the set $K:=K^+\cap
K^-$ consisting of points with bounded total orbits. Each of these sets is
invariant and $K$ is compact. As is in the one dimensional case, recurrence can
occur only on the set $K$. Corresponding to the  set
$J$ in dimension one, we have the sets $J^\pm:=\partial K^\pm$, and the set
$J:=J^+\cap J^-$. Each of these sets is invariant and $J$ is compact.
 A normal families argument shows that there is no
``forward'' instability in the interior of $K^+$ and no ``backward''
instability in the interior of $K^-$.
Thus ``chaotic'' dynamics, that is recurrent dynamics with instability in both
forward and backward time, can occur only on the set $J$.
 
   The techniques that Fatou and Julia used in one dimension are based on
Montel's theory of normal families and do not readily generalize to higher
dimensions. A  different tool appears in the work of Brolin [Br], who made use
of the theory of the logarithmic potential. Potential theory associates to any
compact subset of the plane a measure which is called the harmonic or
 equilibrium
measure, and the ``potential" of this measure which is called the Green
 function.
Brolin showed that for a polynomial map of $\cx{}$ there is an explicit
 dynamical
formula for the Green function. He proceeded to show that the harmonic measure
 of
the Julia set is an invariant measure with interesting dynamical properties.  It
was later observed that potential theory provides alternate proofs of many of the
basic facts of Fatou-Julia theory (see [Si], [T], and [C]).
 
Potential theory in one variable has a natural extension to several complex
variables called pluripotential theory (cf. Klimek [Kl]). In this context the
analogs of the Green function corresponding to the sets $K^+$ and $K^-$ are the
functions $G^+$ and $G^-$. These functions were studied by J.H.\ Hubbard from a
topological viewpoint (see [H] and [HO]). N.~Sibony had the idea of introducing
potential theory into the study of these two-dimensional mappings, he
introduced the two (1,1) currents $\mu^\pm=(2\pi)^{-1}dd^cG^\pm$ and the
measure $\mu=(4\pi)^{-1}(dd^c(G^+\vee G^-))^2$.  Bedford and Sibony established
some properties of $\mu^\pm$ and $\mu$, the results they obtained are contained
in \S3 of [BS1].  (See also [Be].) Further results are contained in 
[BS2--4], [FS]. In the pluri-potential context, $\mu$ is the analogue of the
equilibrium (or harmonic) measure of the set $K$ (and also of $J$). Hubbard and
Papadopol [HP] have shown that a current like $\mu^+$ also arises naturally from
a (non-invertible) holomorphic mapping $f:{\bf P}^n\to{\bf P}^n$.
 
In this paper we combine potential-theoretic methods with tools
from ergodic theory, especially Pesin's theory of non-uniform hyperbolicity.
These tools allows us to describe the geometric structure of the
currents $\mu^\pm$ and to give a geometric description of the relation
between $\mu^\pm$ and $\mu$. The starting point for these results is a
characterization of the measure $\mu$ in terms of entropy which we now
describe.
 
We can associate to each invariant probability measure $\nu$ its
measure-theoretic entropy $h_\nu(f)$. The variational principle states that
the supremum of $h_\nu(f)$ taken over the set of all invariant
probability measures is the
topological entropy, $h_{top}(f)$. A measure $\nu$ for which
$h_\nu(f)=h_{top}(f)$ is called a {\it measure of maximal entropy}.
For polynomial maps in one dimension the topological entropy is
$\log d$ where $d$ is the degree of the polynomial (see [G] and [Lyu1]),
and $\mu$ is the {\it unique} measure of maximal entropy (see [Lyu2] and [Ma]).
In two complex dimensions the topological
entropy is $\log \deg f$ (see [FM] and [S]), and  $h_\mu(f)=\log\deg f$ (see
[BS4]). In \S 3 we prove: {\sl The harmonic measure is the unique measure of
maximal entropy for a polynomial diffeomorphism of $\cx2$} (Theorem
 3.1).
 
For polynomial maps of $\C$, Fatou and Julia used Montel's theorem
to show that expanding periodic points are dense in $J$. This result can also be
proved using potential theory. A key observation in such a
potential-theoretic proof is the fact that the support of
harmonic measure is the set $J$.  For polynomial diffeomorphisms of $\cx2$ the
situation is not so straightforward. If $J^*\subset\cx2$ denotes the
support of $\mu$, then it follows easily that $J^*\subset J$. For polynomial
diffeomorphisms which are hyperbolic, we have shown in [BS1] that $J=J^*$.  But
the question of whether equality holds in general seems to be very difficult.
 
Periodic saddle points are the analogs of expanding periodic points for two
dimensional diffeomorphisms. These are points of period $n$ for which $Df^n$ has
one eigenvalue outside and one eigenvalue inside the  unit circle. It is
relatively easy to show that every saddle orbit is contained in $J$. In \S 9 we
prove the more difficult result: {\sl Every saddle orbit is contained in $J^*$.}
It was shown in [BS4] that the closure of the saddle orbits contains $J^*$.
Combining these results gives: {\sl $J^*$ is the closure of the set of saddle
orbits}. Thus $J^*$ plays a role for polynomial diffeomorphisms of $\cx2$
analogous to the role played by $J$ for polynomial maps of $\C$.
 
Let $p$ be a periodic saddle point. The stable/unstable manfolds of $p$
are defined as
 $$W^{s/u}(p):=\{q\in\cx2:\lim_{n\to\infty}dist(f^{\pm n}q,f^{\pm n}p)=0\}.$$
In \S2 we show that for $\mu$ almost every point $p$, the set $W^{s/u}(p)$ is
conformally equivalent to $\C$ and is a dense subset of $J^\pm$. This result
was obtained independently by Wu in [W].

For distinct periodic saddle points, $p$ and $q$, the intersections of $W^s(p)$
and $W^u(q)$ are called heteroclinic intersections. We show in \S9 that $J^*$
 can
be characterized in terms of heteroclinic intersections. {\sl For any pair of
periodic saddle points $p$ and $q$: $J^*=\overline{W^s(p)\cap W^u(q)}$}. It is
interesting to contrast this description of $J^*$ with a similar
description of $J$ from [BS4]. {\sl For any pair of periodic saddle points:
$J=\overline{W^s(p)}\cap\overline{W^u(q)}.$}
The intersections of $W^s(p)$ and $W^u(p)$ other than $p$ itself are called
homoclinic intersections.  It was observed in [BS4] that the set of periodic
saddle points that create homoclinic intersections is dense in $J^*$. In \S9 we
prove the more delicate result that {\it every} periodic saddle point creates
homoclinic intersections.
 
The harmonic measure $\mu$ and the currents $\mu^\pm$ are related by the
analytic equation $\mu=\mu^+\wedge\mu^-$. This formula does not give much
geometric insight into the relation between these objects. The results on
periodic saddle points and stable manifolds are consequences of a geometric
description of the currents $\mu^\pm$ and the way in which these currents
``intersect'' to give $\mu$.
In order to explain the results of this paper about general polynomial
diffeomorphisms it is useful to recall results from [BS1] about the special case
of uniformly hyperbolic polynomial diffeomorphisms.
 
A polynomial diffeomorphism $f$ is uniformly hyperbolic if there is a hyperbolic
splitting of the tangent bundle over $J$.  Hyperbolicity implies that for every
point $p\in J$ the sets $W^{s/u}$ are immersed submanifolds. In the uniformly
hyperbolic case, the collection of stable manifolds has the following
``laminar'' structure.  At a point $p\in J$, we may let $T^u$ be a small complex
disk transversal to $W^s(p)$.  For points $q\in J$ near $p$, the local stable
manifold $W_\epsilon^s(q)$ will  intersect $T^u$ in a unique point $a\in T^u$.
If we let $A^u\subset T^u$ denote the set of such intersections, then we may
parametrize the local stable manifolds by $a\in A^u$, and locally $J^+$ is
topologically equivalent to the product of $A^u$ and a disk.  Given two such
transversals $T_1$ and $T_2$ and corresponding sets $A_j\subset T_j$, $j=1,2$,
there is a (continuous) {\it holonomy} map $\chi:A_1\to A_2$, defined by
following a stable disk from its intersection point $a_1\in A_1$ to the point
$a_2\in A_2$ where it intersects  $T_2$.  This gives a homeomorphism between the
intersections with nearby transversals. In [BS1] we showed that the holonomy map
preserves the slice measures $\mu^+|_{T_j}$.
 
There is a corresponding theory, due to Pesin, of (non-uniform) hyperbolicity
with respect to a measure $\nu$. An (ergodic) measure $\nu$ is said to
be hyperbolic if no Lyapunov exponent is zero. (See \S2  for the relevant
definitions.) The theory of Pesin for a hyperbolic
measure $\nu$ implies that for $\nu$ almost every point $p$ the sets
$W^{s/u}$ are immersed  submanifolds. 	It is shown in [BS4] that the measure
$\mu$ is ergodic and hyperbolic.
In the case of a hyperbolic measure,  we may define a holonomy map on a
compact set of positive measure (but not necessarily everywhere).  In \S4 we
show: {\sl
The holonomy map preserves the slice measures $\mu^+|_{T_j}$}.
 
In the uniformly
hyperbolic case, we may take a similar transversal $T^s\subset W^s(p)$, and
we may parametrize the local unstable manifolds by $A^s\subset T^s$.  It
follows that a neighborhood in $J$ is homeomorphic to $A^s\times A^u$.  In the
case of a hyperbolic
measure, it is possible to find product sets with positive measure, which we
call {\it Pesin boxes} and denote again as $A^s\times A^u$. The measure $\mu$
induces conditional measures on each stable slice. As a byproduct of the
characterization of $\mu$ as the unique measure of maximal entropy in \S3 we
show: {\sl The conditional measures on the stable/unstable slices are given by
$\mu^{+/-}$.}   As a consequence of the holonomy invariance of $\mu^\pm$ and the
identification with the conditional measures, we obtain in \S4 the result: {\sl
$\mu$ restricted to a Pesin box is a product measure.} This allows us to invoke
results of Ornstein and  Weiss which imply that {\sl $\mu$ is Bernoulli}. This
is the strongest mixing property that a measure can possess.
 
Now let us pass from the analysis of the slice measures of  $\mu^\pm$ to the
currents themselves.  A closed manifold $M$ defines a current of integration,
denoted by $[M]$.  (See \S5 for a general discussion of currents.)  In the
uniformly hyperbolic case, the laminar structure of $\cW^{s/u}$ passes over to a
laminar structure for $\mu^\pm$.  That is, at a point $p\in J$, we may choose an
open set $U$ and a transversal $T^u$ such that for each $a\in A$, the local
stable manifold $D^s(a)$ is a closed submanifold of $U$, and the restriction of
$\mu^+$ to $U$ is given by
\footnote{*}{Throughout this paper we use the following notation for
integration.  If $\lambda$ is a measure on $A$, and if $f$ is an
integrable function on $A$ with values in the space of currents, then we write
the integral as $\int_{a\in A}\lambda(a)\,f(a)$.}
$$\mu^+\contract
U=\int\lambda^u(a)\,[D^s(a)],\eqno(\ddag)$$ which is a direct integral of
currents of integration with respect to $\lambda^u$, which is the measure
obtained by restricting $\mu^+$ to $T^u$.  A current of the form (\ddag) is
called {\it uniformly laminar} if the manifolds $D^s(a)$ are pairwise disjoint.
With the family $\cW^{s/u}$ given by Pesin theory, there is no uniformity to the
size of the manifolds, i.e.\ we cannot choose $U$ such that for $M\in\cW^{s/u}$
every component of $M\cap U$ is closed in $U$.  In \S6 we define the more
 general
class of {\it laminar} currents and show that a laminar current $T$ is given as
a countable sum, $T=\sum T_j$, where the $T_j$'s have disjoint carriers, and
$T_j$ is uniformly laminar on some (possibly small) open set $U_j$.  The closed,
laminar currents give a natural generalization of the current of integration and
seem to be an interesting class in their own right.  In \S7, it is shown that
$\mu^\pm$ is laminar.  In \S8, it is shown that $\mu^+$ contains uniformly
laminar ``pieces" whose structure is induced by the Pesin boxes.  This allows us
to show that the wedge product that defines the measure $\mu$ is in fact given
by an intersection product of stable and unstable manifolds. This yields further
structure for the currents $\mu^\pm$.
 
Much attention has been paid to polynomial diffeomorphisms with real
 coefficients.
In this case the real subspace $\RR^2\subset\CC^2$ is invariant and we may let
$f_\R$ denote the
 restriction of $f$ to
$\RR^2$. The (real) H\'enon map is a much studied example with $\deg=2$. In
 contrast to the complex
case, where the topological entropy is $\log d$, the topological entropy of
$f_\R:\RR^2\to\RR^2$ can be any real number in
 the interval
$[0,\log d]$ (see [FM] and [Mi]).  In \S10 we give several equivalent criteria
 for the
entropy of $f_\R$ to be equal to $\log d$.  One of these is that $K\subset
 \RR^2$, that
is to say that every complex bounded orbit is actually real. A second criterion
 is that
every periodic point of $f$ is in $\RR^2$.  A third is that: {\sl For any
 hyperbolic
point $p$, all intersection points $W^s(p)\cap W^u(p)$ lie inside $\RR^2$}.
These results may be used to show that, when topological entropy is maximal, the
loss of a single periodic point or homoclinic intersection forces a decrease in
the the topological entropy.

This paper is divided into different parts, according to the methods that
predominate.  In \S\S2--4 the principal tools are Smooth Ergodic Theory,
especially Pesin's Theory.
 In \S\S5--7, the primary tools are the theory of currents and
the Ahlfors Covering Theorem. These sections do not use Ergodic Theory.
Finally, these
methods are combined in \S\S8--9.
 
The specific contents are as follows.  \S2 gives a summary of the part of
Smooth Ergodic Theory that we will use.  At the end of \S2 it is shown that,
for $\mu$ almost every point $p$, the stable manifold of $p$ is dense in $J^+$ and
conformally equivalent to $\cx{}$.  In \S3 the conditional measures are shown to
be induced by the current $\mu^+$.  (This permits estimates on the Hausdorff
dimension and Lyapunov exponent at the end of the section.)  Then it is shown
that $\mu$ is the unique measure of maximal entropy.  The holonomy map is
discussed in \S4, and it is shown that the holonomy of the Pesin stable manifolds
preserves the restriction measures of $\mu^+$.  Finally, it is shown that $\mu$
has a local product structure.  In \S5 we summarize the main ideas and
definitions that we use from the theory of currents.  Laminar currents are
defined in \S6, and the basic structure is developed.  In \S7 it is shown that
$\mu^\pm$ are laminar currents.  In \S8 we show that the laminar structure of
$\mu^+$ coincides with the structure induced by the Pesin manifolds and the
conditional measures.  And in \S9 we apply the previous work to the
 study
of saddle points.   Real H\'enon mappings are discussed in \S10, and
several (equivalent) criteria are given for $f$ to be essentially real. \S11 is
an appendix which outlines an alternative sequence in which the results of
this paper can be obtained.  This alternate approach starts with results of
Pesin theory and then proceeds to the theory of currents.  The main difference
is that the use of the methods of entropy theory is delayed until the
end.
 
\bigskip
 
\section 2.  Preliminaries from Ergodic Theory

 {\bf 2.1. Measurable Partitions
and Conditional Measures}\medskip\noindent
The technique of measurable partitions developed by Rokhlin [Ro1] is a
powerful tool in measure theory. Somehow it is not widely known beyond
ergodic theory. So, we will spend some time to define the main concepts and
to establish notation.

Let $J$ be a compact metric space, and let $\nu$ be a  probability Borel measure
 on $J$.
A {\it partition} $\xi=\bigcup\xi_\alpha$ of $J$ is a decomposition of $J$
into disjoint, measurable subsets.  The element of the partition
containing $x$ will be denoted by $\xi(x)$, and will be called the {\it fiber}
through $x$.  Note that all fibers can have zero measure. For example we can
consider a partition $\varepsilon$ into single points.
Two partitions are considered to be equivalent
if they coincide on a subset $J^\prime$ of full measure.
 
Each measurable function $\phi$ generates  a partition whose fibers are level
sets of $\phi$. Such partitions are called {\it measurable}. Any countable
partition is measurable.
 An orbit partition of an irrational rotation of the circle (with Lebesgue
measure) gives an example of non-measurable partition.  More generally, one can
consider an orbit partition of any ergodic transformation; see the discussion
below.
 
The basic property of measurable partitions is for any measure $\nu$ there is a
{\it family of  conditional measures} $\nu(\cdot|\xi(x))$ on the fibers. This
family is uniquely determined by the following properties:
 
\smallskip
\item{(i)}  Each $\nu(\cdot|\xi(x))$  is a probability measure on $\xi(x)$;
\item{(ii)} For any integrable function $\phi$, the function
$${\phi_\xi(x)=\int\phi(y)\, \nu(y|\xi(x))}$$
along the fibers is measurable and integrable, and
 
\item{(iii)} $$\int{\phi_\xi(x)} \, \nu(x) =\int\phi\, \nu.$$
 
\medskip
\give Remark. The above averaging of $\phi$ over the conditional measures
is equivalent to taking of the conditional expectation of $\phi$
with respect to the $\sigma$-algebra generated by $\xi$.\medskip
 
By ``countable" set we will mean ``at most countable".
If we have a countable family of measurable partitions $\xi_i$ then we can
construct a partition $\bigvee\xi_i$ by intersecting fibers of $\xi_i$, i.e.
$$(\bigvee\xi_i)\,(x)=\bigcap\,\xi_i(x).$$
One can check that this construction leads to a measurable partition.
 
Finally, let us mention that
for an arbitrary (non-measurable) partition $\eta$ there exists
its {\it measurable envelope},
i.e. the finest measurable partition which is coarser than $\eta$.

\bigskip
\noindent
{\bf 2.2. Elements of Entropy Theory}\medskip
 \noindent
The reader can see [Ro2] or [CFS] for the background in entropy theory.
Our exposition will be adapted to our goals (in particular, it will not be
as general as possible).
 
Entropy of a countable (mod 0) measurable partition $\xi=\{\xi_i\}$ is
defined as
$$H_\nu(\xi):=-\sum\nu(\xi_i)\log\nu(\xi_i)
= \int\log {1\over \nu(\xi(x))}\nu(x)\eqno (2.1)$$
(it can be infinite). If
the partition is not countable then its entropy is infinite by definition.
 
If we have two measurable partitions $\xi$ and $\eta$ then we can
  restrict $\xi$
on the fibers of $\eta$, thus we can calculate the entropy of $\xi$ with
respect to $\eta$ in terms of the conditional measures as
$H_{\nu(\cdot|\eta(x))}(\xi|\eta(x))$.  We then define the {\it conditional
entropy} by averaging this with respect to $\nu$:
 
$$H_\nu(\xi|\eta):=\int\, H_{\nu(\cdot|\eta(x))}(\xi|\eta(x))\, \nu(x).$$
 
Let us consider now a homeomorphism $f : J\rightarrow J$ preserving  a measure
$\nu$. Then it naturally acts on the space of measurable partitions
$\xi\mapsto f\xi$, where $(f\xi)(x)=f(\xi(f^{-1}x))$. A partition $\xi$ is
 called
$f$-{\it invariant} if $f\xi$ is a refinement of $\xi$. A partition $\xi$ is
called a {\it generator} if
$$\bigvee_{n=-\infty}^\infty  f^n\xi= \varepsilon.$$

Given a partition $\xi$, consider the $f^{-1}$-invariant partition
$\xi^u=\bigvee_{n=0}^\infty f^n\xi$. Let us call the fibers of this partition
$\xi$-{\it unstable} fibers. We can define the Jacobian $J^u f$ of $f$
in the ``$\xi$-unstable direction" as the Radon-Nikodym derivative of $f$ with
respect to conditional measures:
$$J^u f(x)= {df^* \nu(\cdot|\xi^u(fx))\over d\nu(\cdot|\xi^u(x))}.$$

Since $\nu$ is invariant, $J^u f$ is constant on the fibers of $f^{-1}\xi^u$,
and hence
$$J^u f(x)={1\over p(x) }\eqno (2.2)$$
where $p(x)=\nu(f^{-1}\xi^u (fx)| \xi^u(x))$.
Now define entropy of $f$ with respect to $\xi$ as
$$h_\nu (f,\xi)= H_\nu (f^{-1} \xi^u | \xi^u) =
-\int\log p(x)\nu(x) =
\int\log J^u f(x) \nu(x)\eqno (2.3)$$
(the middle equality follows from (2.1)). So, from the dynamical
point of view entropy of a transformation with respect to a partition
is just the logarithm of the geometric average of the Jacobian of $f$ in the
$\xi$-unstable direction.
 
Finally, the entropy of $\nu$ with respect to $f$ is defined as
$$ h_\nu(f)=\sup_\xi h_\nu(f,\xi)$$
where supremum is taken over all measurable partitions $\xi$. Actually, one
can take the supremum over finite partitions only. Moreover, it is
enough to evaluate entropy of any generator with finite entropy:
\proclaim Proposition 2.1. If $\xi$ is a generator with finite entropy
    then $h_\nu(f)=h(f,\xi)$.
 
In conclusion let us discuss the ergodic decomposition of the transformation
 $f$.
Let us consider the {\it orbit partition}
 $O$ of $f$ (whose fibers are orbits of $f$).
The transformation $f$ is ergodic if the measurable envelope of
$O$ is a trivial partition (whose only fiber is the whole space).
 
In general, let us consider the measurable envelope $E$ of $O$.
The fibers of $\varepsilon$ supplied with conditional measures are
called {\it ergodic components} of
$\mu$.   Note that all ergodic
components may have zero measure (consider the identity transformation). However
it makes sense to consider the whole space of these components mod 0.
Restricting $f$ onto ergodic components and then taking them together we obtain
a representation of $f$ as a ``direct integral" of ergodic transformations.
It follows from (2.3) that
$$  h(f)=\int h(f|E(x))\nu(x).\eqno (2.4)$$
This formula gives a method for reducing entropy questions for arbitrary
measures to the case of ergodic measures.  If $f$ is ergodic then it has a
 finite
generator by the Krieger Theorem [Kr]. So, in the ergodic case we can always
compute entropy according to  Proposition 2.1.
 
For ergodic $\nu$ let us say that a point $x$ is $\nu$-equidistributed if for
any continuous function $\phi$
$$\lim_{n\to\infty}{1\over n}\sum_{k=0}^{n-1} \phi (f^k x)=\int  \phi \nu .$$
By the Birkhoff Ergodic Theorem $\nu$ almost every point is
$\nu$-equidistributed.

\bigskip
\noindent
{\bf 2.3. Measures of maximal entropy}\medskip\noindent
For the material of this section we refer to Bowen's book [Bo].
We will not define the topological entropy of $f$, but a basic property is
given by the so-called Variational Principle, which asserts that the
topological entropy $h(f)$ is given as
$$h(f)=\sup h_\nu(f)\eqno (2.5)$$
where $\nu$ runs over all probability
Borel measures invariant with respect to $f$.
 
A measure $\mu$ is called a {\it measure of maximal entropy} if
$h_\mu(f)=h(f)$. This measure does not neccessarily exist, but if it does,
then by (2.4) all its ergodic components are measures of maximal entropy as
well. Hence, existence/uniqueness of a measure of maximal
entropy are equivalent to the existence/uniqueness of an {\sl ergodic}
measure of maximal entropy.
 
The problem of uniqueness of the measure of maximal entropy is not
handled yet in a general setting. The status of the existence problem
is much better:
 
\proclaim Newhouse Theorem [Ne]. If $f: M\rightarrow M$ is a
 $C^\infty$
diffeomorphism of a compact $C^\infty$-manifold then $f$ has a measure of
maximal entropy.

\bigbreak
\noindent
{\bf 2.4. Stable and unstable manifolds}\medskip

Some basic references for the material in this section are [P1], [FHY], [R2] and
[PS].
Let $M$ be a  Riemannian $C^2$-manifold,
$f: M\rightarrow M$ be a $C^2$-diffeomorphism,
 $J$ be  an invariant compact subset  of $M$.
Let $\nu$ be an invariant {\it ergodic}  measure of $f$ supported on $J$.
As usual $T_x M$ denotes the tangent space at $x$.   A measurable function
$r(x)$ is called $\epsilon$ {\it slowly varying} if
$$(1+\epsilon)^{-1} r(x) < r(fx) < (1+\epsilon) r(x).$$
 
\proclaim Oseledec Theorem. There exist finitely many distinct real numbers
$\chi_i,\; i=1,...s$ called
characteristic exponents, an invariant set ${\cal R}$ of full measure,
and $s$ invariant measurable  distributions
 $E_i(x)\subset T_x M,\; x\in {\cal R}$,
such that\smallskip
\item {(i)} $T_x M=\oplus E_i (x) $;\smallskip
\item {(ii)} For any nonzero $v\in E_i(x)$, $$\lim_{n\to\pm\infty}{1\over n}\log
\| Df^n (x) v\| =\chi_i.$$
\item{(iii)}  For $i\ne j$ and $\epsilon>0$ there is an $\epsilon$ slowly
varying function $s_{ij}(x)>0$ which is less than the angle between $E_i(x)$
and $E_j(x)$.
 
The points of the set ${\cal R}$ are called {\it regular}. We can also assume
that ${\cal R}$ consists of $\nu$-equidistributed points.

Let us state now the Pesin Theorem which says that the above distributions are
integrable.
Denote by $B(x,r)$ a ball of radius $r$ centered at $x$,
and $B^{s/u}(x,r)=E^{s/u}(x,r)\cap B(x,r)$. Now let us
define {\it stable} and {\it stable-center} distributions
 
$$E^s(x)=\oplus_{\chi_i<0} E_i,\quad E^{sc}(x)=\oplus_{\chi_i\leq 0} E_i.$$
Similarly one can define the unstable and unstable-center distributions
 $E^u(x)$ and $E^{uc} (x)$.

\proclaim Pesin Theorem [P1]. Let $dim \,E^s>0$. Then for any
$\epsilon>0$ there are $\epsilon$-slowly varying positive functions
$C(x)=C_\epsilon(x)$ and
$r(x)=r_\epsilon(x)$ on ${\cal R}$, and a
family $W_{loc}^s(x),\; x\in {\cal R}$ of smooth  manifolds satisfying the
following properties:\smallskip
\item {(i)} $W_{loc}^s (x)$ is a graph of a function
 $B^s (x, r(x))\rightarrow E^{uc} (x)$ tangent to $E^s (x)$ at $x$;\smallskip
\item {(ii)} For any $y\in W^s_{loc} (x)$  and $n=1,2,\dots$
$$C(x)^{-1} \exp (-(\lambda+\epsilon)n) \leq \ {\rm dist} (f^n x,
f^ny) \leq C(x) \exp (-(\lambda-\epsilon)n);$$ \item
{(iii)} The $f$ underflows the manifolds $W^s_{loc}(x)$:
   $f W^s_{loc}(x)\subset W^s_{loc}(fx).$
 
The manifolds $W^s_{loc}(x)$ are called {\it local stable manifolds}.
For $r\leq r(x)$ let $W^s_r (x)$ be a part of $W^s_{loc} (x)$ lying over
$B^s (x,r)$.
In order to obtain the Theorem on {\it local unstable manifolds}
$W^u_{loc}(x)$ we just interchange the roles of $f$ and $f^{-1}$.
 
Let us indicate one immediate consequence of this result.
 
\proclaim Proposition 2.3.  If the measure $\nu$ is not supported on a periodic
orbit then all characteristic exponents cannot be negative
(positive).
 
\give Proof.  Otherwise $W^s_{loc}(x)=B(x, r(x))$. Since
$\nu$-equidistributed points
$x\in {\cal R}$ are recurrent, we can find a moment $n>0$ such that
$f^n$ maps $B(x,r)$ into itself uniformly contracting it. It follows that
$x$ is periodic, and $\mu$ is supported on its orbit. \qed
 
The family of local unstable manifolds does not form a partition.
The following statement supplies us with a $f^{-1}$-invariant measurable
 partition (called a {\it Pesin partition}) subordinate
to the family of manifolds.
 
\proclaim Theorem  2.4 (see [P2], [LS]).
There is a measurable $f^{-1}$-invariant generator
 $\xi^u$  whose fibers are open subsets of the local unstable manifolds, and
such that $$h_\nu(f)=h_\nu(f,\xi^u).$$

\give Remark. When $f$ has no zero characteristic exponents then the Pesin
partition $\xi^u$ is a $\xi$-unstable partition for some partition $\xi$
with finite entropy [LY]. So, in this case the above entropy formula follows
from Proposition 2.1.
 
Let us now define the {\it global unstable manifold} $W^u(x)$ at $x$ as the set
of points $y$ whose backward orbits are asymptotic to the backward orbit of $x$.
Clearly $f W^u(x)= W^u (fx).$
One can prove that for $x\in {\cal R}$
$$W^u(x)=\bigcup f^{n} W^u (f^{-n} x).$$
This implies the following two consequences:

\item{(i)} The backward orbits $y\in W^u(x)$ are  exponentially
asymptotic to the orbit of $x$.
 
\item{(ii)} The set $W^u (x)$  is an immersed Euclidean space.
 
The global unstable manifolds form the partition of the measure space
$(J, \nu)$ which we will call the {\it global unstable partition}.
This partition is in general not measurable.
 
A partition $\tau$ is called {\it hyperfinite} if there is
 a sequence of measurable partitions $\xi_i$ such that
$$\xi_1(x)\subset\xi_2(x)\subset...,\; {\rm and}\; \tau(x)=\bigcup\xi_i(x).$$
Let us call a measure defined up to a scalar factor a {\it projective
measure class}.
On a fiber of a hyperfinite partition one can define a conditional
projective measure class ${\bar\nu}(\cdot|\tau(x))$ as the class of the
measure:
$$\nu(\cdot|\tau(x))=\lim_{i\to\infty}
{\nu(\cdot|\xi_i(x))\over\nu(\xi_1(x)|\xi_i(x))}.$$ 
Thus for any other
sequence $\xi_n^\prime$ which generates $\tau$, the measure
$\nu^\prime(\cdot|\tau(x))$ obtained in this way will be a multiple of
the measure above by a constant depending only on $x$.   In fact, for any
measurable partition $\eta$ subordinate to $\tau$, the conditional measures on
the fibers of $\eta$ are just the normalized projective measure classes of
$\tau$.

\proclaim Proposition 2.5. The global unstable  partition is hyperfinite.
 
\give Proof.  Take a Pesin partition $\xi^u$, and represent the global
unstable partition as the limit of measurable partitions $f^{-n}\xi^u$. \qed
 
It is evident that the preceding discussion may be applied equally well to the
stable direction instead of the unstable one.
 
\bigskip
\noindent
{\bf 2.5. Relations between entropy and characteristic exponents}\medskip
 \noindent
The following inequality was discovered by Margulis in
the case of an absolutely continuous measure. It was later generalized
by Ruelle [R1]:
 
\proclaim Margulis-Ruelle inequality.
$$ h_\nu (f)\leq \sum_{\chi_i>0} \chi_i \dim E_i, $$
and a corresponding inequality holds with the sum of
negative characteristic exponents.
 
\proclaim Corollary. If $h_\nu(f)>0$, and if $f$ has at most two characteristic
exponents, then one of these exponents is negative, and another is positive.
 
In such a situation we will denote the negative and positive exponents $\chi^s$
and $\chi^u$ correspondingly.
 
 More recently a number of remarkable
relations between entropy, characteristic exponents and Hausdorff dimension
have been discovered (see Pesin [P1] and Ledrappier-Young [LY] and the
references there.)
The {\it Hausdorff dimension} of a measure $\nu$, written $HD(\nu)$, is defined
as the infimum of the Hausdorff dimension of $X$, for all Borel subsets $X$
with full $\nu$ measure.  Clearly, the Hausdorff dimension depends on the
measure class only.
 
\proclaim Lai-Sang Young's Formula [Yg].  Assume that $f$ has only one
characteristic exponent $\chi^s<0$.  Then for $\nu$~a.e.\ $x$,
$$HD(\bar\nu(\cdot|W^s(x)))={h_\nu(f)\over|\chi^s|}.$$
 
\bigskip
\noindent
{\bf 2.6. Complex analytic case}\medskip
\noindent 
Let $M$ be a Hermitian complex analytic manifold and let $f$ be
analytic. Then $E_i(x)$ are complex subspaces in $T_x M$, and all local
manifolds are complex analytic.
 
Assume now that dim$_{\bf C} M=2$, and and $\nu$ be any invariant probability
measure with two non-zero characteristic exponents of opposite signs,
$\chi^s<0$ and $\chi^u>0$.
(In particular, this will be the case if $h_\nu (f)>0$,
see  the Corollary of the Margulis-Ruelle inequality).
Hence dim$_{\bf C}E^{s/u}=1$, the global stable/unstable manifolds are
regular complex curves. The following statement says that almost all of them are
parabolic.
 
\proclaim Proposition 2.6.
The stable and unstable manifolds $W^u(x)$ and $W^s(x)$ are
conformally equivalent to the complex plane for $\nu$ a.e.~$x$.
 
\give Remark.  It is possible to prove Proposition 2.6 along the
lines of the proof of Theorem~5.4 of [BS1]. That is, {\sl for $x\in{\cal R}$,
 $W^u(x)$
contains a sequence of disks $D_1\subset D_2\subset\dots$ such that the modulus
of $D_{j+1}-D_j$ is bounded below.}  From this, it follows that $W^u(x)$ is
equivalent to $\cx{}$.  To sketch this argument, we note that
$W^u_{r(f^{-k}x)}(f^{-k}x)$ is a graph over a disk of radius $r(f^{-k}x)\ge
c(1+\epsilon)^{-k}r(x)$.  On the other hand, the derivative of $f^{-k}$ on
$W^u_{r(x)}(x)$ is approximately $e^{-n\chi^u}$.  For a small disk $D_j$
containing $x$ inside $W^u(x)$, we may choose $n$ sufficiently large that the
modulus of the annulus $W^u_{r(f^{-n}x)}(f^{-n}x)-f^{-n}D_j$ is at least 2.
 Then
we may let $D_{j+1}=f^nW_{r(f^{-n}x)}(f^{-n}x)$.

\give Remark. The proof we give below uses a technique that will also be
used in \S3. Two measurable functions $\alpha$ and $\beta$ are called {\it
cohomologous} if there is a measurable function $\omega$ such that the following
cohomology equation
$$ \alpha(x)-\beta(x)= \omega(fx)-\omega(x)$$
is satisfied $\nu$-almost everywhere.
 
Usually, the cohomology equation comes up when we calculate the logarithm
of the Jacobian (or norm) of $f$ with respect to two equivalent measures
(metrics). The following statement will be useful on several occasions.
 
\proclaim Lemma 2.7. Let $\alpha$ be a measurable function bounded from below.
If $\alpha$ is cohomologous to 0 then
$$\int\alpha \,\nu = 0.$$
 
This is trivial if $\omega$ is integrable. Otherwise, the proof is based
upon the Birkhoff ergodic theorem (see, e.g., [LS, Proposition 2.2].)
 
 \medskip

\give Proof of Proposition 2.6. We consider the unstable manifolds
$W^u(x)$.
 Let $F_{W^u(x)}$ denote the
Kobayashi metric on $W^u(x)$.  This
metric depends in a lower semicontinuous manner on $x$ if $W^u(x)$ depends
continuously on $x$.  And since we may find compact subsets of $J$ of
measure arbitrarily close to 1 on which $W^u(x)$ depends continuously on $x$,
the correspondence $x\mapsto F_{W^u(x)}$ is measurable.
$W^s(x)$ is conformally equivalent to either a plane
or a disk, depending on whether $F_{W^u(x)}(x,E^u)=0$ or not.  By ergodicity,
 the
type of $W^u (x)$ is the same for almost all $x$. We assume that it is
hyperbolic and derive a contradiction.
 
Let $\alpha (x)=\log \|Df(x)|_{E^u}\|$, where the norm is taken
with respect to the Hermitian metric on $M$.  Similarly, for $x\in{\cal R}$ we
define the function $\beta (x)=\log |Df(x)|_{E^u}|$, where
$|Df(x)|_{E^u}|$ denotes the norm taken with respect to the Kobayashi
metric.  If we let $\rho(x)$ denote the ratio of the Hermitian to the
Kobayashi metrics in the unstable direction, then $\alpha$
and $\beta$ are cohomologous in the sense that
$$\alpha(x)-\beta(x)=
\log\rho(fx)-\log\rho(x).$$
But $f$ is an isomorphism between $W^u(x)$ and $W^u(fx)$ and hence preserves the
Kobayashi metric. So $\beta(x)=0$ almost everywhere, and  $\alpha$ is
cohomologous to 0. By Lemma 2.7,
$$\chi^u\equiv\int \alpha \,\nu = 0$$
contradicting  the assumption that $\chi^u>0$. \qed
 
 Assume now that $f$ is a polynomial automorphism of ${\bf C}^2$, and $J^+$,
$J^-$,  $J=J^+\cap J^-$ etc., be the sets introduced in \S 1. Consider also the
currents $\mu^+$ and $\mu^-$.
We can ``slice'' the current $\mu^+$ with any
complex one dimensional variety $W$ (see \S5).
The result can be interpreted as a measure
on $W$ which we denote by $\mu^+|_W$.
In particular, we can consider the measure $\mu^+|_{W^u(x)}$
on the unstable leaf containing $x$. We will call it an {\it unstable slice} of
$\mu^+$.
 
\proclaim Lemma 2.8. Any $\nu$ regular point $x\in {\cal R}$ belongs
to the support of $\mu^+|_{W^u(x)}$.
 
\give Proof.  Let $\Delta$ be a disk with $x\in\Delta\subset W^u(x)$. If
$\mu^+(\Delta)=0$ then $G^+$ is harmonic in $\Delta$. Since the orbit of $x$ is
bounded $G^+(x)=0$ and by the minimum principle $G^+$ is zero in $\Delta$. Thus
$\Delta\subset K^+$. It follows that $f^n(\Delta)\subset K^+$ so $f^n(\Delta)$
is uniformly bounded for all $n$. By the Schwartz lemma $||D(f^n|_\Delta)||\le
C$ but this contradicts the fact that at a regular point the Lyapunov exponent
is positive in the unstable direction. \qed
 
\proclaim Proposition 2.9. For $\nu$ a.e.~$x$, $W^u(x)$ is a dense subset
of $J^-$, and $W^s(x)$ is a dense subset of $J^+$.
 
\give Proof.  For $x\in {\cal R}$ it is evident that $W^s(x)\subset K^+$.
We show that $W^s(x)\subset J^+$.  Let us suppose that $y\in
W^s(p)\cap int K^+$. Since the iterates $f^n$ for $n\ge 1$ form a normal family,
it follows that
$||Df^n_y||$ is bounded. The fact that $y\in W^s(x)$ implies that
 $d(f^n(x),f^n(y))\le
Cr^n$ for $r<1$. This in turn implies that $||Df_{f^n(x)}-Df_{f^n(y)}||\le
C'\rho^n$ with some $\rho < 1$. It
follows from [R, Theorem 4.1] that the asymptotic behavior of $Df^n_x$ and
 $Df^n_y$ is
the same. In particular
$$\lim_{n\to\infty}{1\over n}\log||Df^n_y||=\lim_{n\to\infty}{1\over
n}\log||Df^n_x||=\chi^u$$
We see that $||Df^n_y||$ is therefore not bounded. This
completes the proof that $W^s(x)\subset J^+$.
 
We now show that that $W^s(x)$ is dense in $J^+$.
Let $D^s$ denote a disk inside $W^s(x)$ with $x\in D^s$ and such that
$\mu^-|_{W^s(x)}(\partial D^s)=0$.  By [BS3], $d^{-n}[f^{-n}D^s]$ converges to
$c\mu^+$ as $n\to\infty$, with
$$c=\mu^-|_{W^s(x)}( D^s).$$
By Lemma 2.8 (in the ``stable" setting) $c>0$.

Now let $U$ be an open set with $U\cap J^+\ne\emptyset$.  Thus $\mu^-\contract
U\ne0$, and  so $(f^{-n}D^s)\cap U\ne\emptyset$ for all $n$ greater than some
large $N$.  It follows that
$$f^{-n}W^s(x)\cap U=W^s(f^{-n}x)\cap U\ne\emptyset.$$
Let $S_N=\{x\in{\cal R}:W^s(f^{-n}x)\cap U\ne\emptyset{\rm\ for\ }n\ge N\}$.
Clearly, $S_1\subset S_2\subset\dots$, and by the previous remark $\bigcup
S_N={\cal R}$.  Further, $fS_N=S_{N+1}$, and since $f$ is ergodic, each $S_N$
 has
measure 0 or 1.  Thus it follows that $S_0$ has full measure, which completes
the proof. \qed
\give Remark.  If $p$ is a (periodic) saddle point, then the average $\nu$ of
point masses over the orbit of $p$ is a hyperbolic measure.  Thus Theorem 1 of
[BS3] is a consequence of Proposition 2.9.

\section 3. The unique measure of maximal entropy
 
 The goal of this section is to prove the following  Uniqueness Theorem.  Our
approach is reminiscent of the proof of the uniqueness theorem for rational
endomorphisms of $\cx{}$ given in [EL] and also of Ledrappier's proof of the
``Variational Principle" for absolutely continuous invariant measures [Le].  An
important consequence of the proof is that the conditional projective measure
class of $\mu$ on the unstable foliation is induced by the current $\mu^+$. At
the end of the section we will derive estimates of Hausdorff dimension
and characteristic exponents.
 
\proclaim Theorem 3.1.
The measure $\mu$ is the unique measure of maximal entropy.
 
Note that this result gives a characterization $\mu$ in terms of topological
dynamics which makes no reference to potential theory.
 
By Theorem 2.2 and the comment preceeding it, there is an
an ergodic measure $\nu$ of maximal entropy, $h_\nu(f)=\log d$.
We are going to show  that $\nu=\mu$ which yields Theorem 3.1.
In fact, this gives an alternative proof that
$\mu$ is a measure of maximal entropy, which was originally proved in [BS4]
 (yet another approach is outlined in \S11).
 
By the Corollary of the Margulis-Ruelle inequality, $\nu$ has two
non-zero characteristic exponents of opposite signs, $\chi^s < 0$ and
$\chi^u > 0$.  So, we can consider the complex one dimensional unstable
foliation and the projective measure class $\bar\nu(\cdot|W^u(x))$ on
its leaves (see Section 2.4). On the other hand, we can consider
measures $\mu^+|_{W^u(x)}$ induced by the current $\mu^+$ (see Section 2.6). 
  For an
open set $B\subset W^u(x)$, we will use the notation
$\mu^+(B):=\mu^+|_{W^u(x)}(B)$.

\proclaim Proposition 3.2. If $\nu$ is a measure of maximal entropy then
for $\nu$ almost all $x$, the conditional projective measure class
$\bar\nu(\cdot|W^u(x))$  is induced by the current  $\mu^+.$
 
\give Proof of Proposition 3.2.
The Jacobian $J^u_{\mu^+}$ of $f$ with respect to the family
of unstable slices of  $\mu^+$ is
equal to $\log d$ since for any $B\subset W^u(x)$ we have
$$
\mu^+(fB)=f^*\mu^+(B)=d
$$
Let $\xi^u$ be the unstable Pesin partition for $\nu$.
By Lemma 2.8  $\rho(x)\equiv\mu^+(\xi^u(x))>0$.
So, we can normalize the above  family of measures in order to get probability
measures on the Pesin pieces:
$$\eta(B|\xi^u(x))={\mu^+(B)}\big/\rho(x).$$
Then the Jacobian $J^u_\eta$ is multiplicatively cohomologous to the
Jacobian $J^u_{\mu^+}$, that is:
$$J^u_\eta(x)=d{{\rho(x)}\over{\rho(fx)}}.$$
So,  $\log J^u_\eta(x)$ is (additively) cohomologous to $\log d$:
 $$\log J^u_\eta(x)-\log d = \log\rho(x) - \log\rho(fx),$$
 
This formula and the following property of the function $\rho(x)$ imply that
$\log J^u_\eta(x)$ is positive.
 
\give Claim. $\rho(f(x))\le d\rho(x)$.
 
By the increasing property of $\xi^u$, we have
$(f^{-1}\xi^u)(x)\subset\xi^u(x)$. So $$\mu^+((f^{-1}\xi^u)(x))\le
\mu^+(\xi^u(x)).$$ On the other hand $f((f^{-1}\xi^u)(x))=\xi^u(f(x))$. So
$$\mu^+(\xi^u(f(x)))=d\mu^+((f^{-1}\xi^u)(x)).$$ Thus
$\mu^+(\xi^u(f(x)))\le d\mu^+(\xi^u(x))$ as was to be shown.
 
Hence $\log J^u_\eta\geq 0$, and  Lemma 2.7 yields
$$\int\log J^u_\eta(x) \nu(x)=\log d,$$ or
$$-\int\log q(x) \nu(x)=\log d \eqno{(3.1)}$$
where $q(x)$ is the $\eta$-measure of $(f^{-1}\xi)(x)$.
 
On the other hand, set $p(x)=\nu(f^{-1}\xi^u (fx)| \xi^u(x))$.
Then by (2.3) and Theorem 2.4
$$-\int\log p(x) \nu(x) = h_\nu (f,\xi^u) = h_\nu(f) = \log d . \eqno (3.2)$$
 
From (3.1) and (3.2) we
conclude
$$\int\log{{q(x)}\over{p(x)}}\nu(x)=0.$$
But
$$\int{{q(x)}\over{p(x)}}\nu(x)=\int_{X}\left(\sum_{\xi^u(y)/f^{-1}\xi^u}
  {{q(z)}\over{p(z)}}p(z)\right)\nu(y)=1.$$
 By concavity of $\log$, we get $q(x)=p(x)$ almost everywhere. Thus
conditional measures coincide on the partition $f^{-1}\xi^u | \xi^u(y)$
for $\nu$ almost all $y$. The
same argument applied to $f^n$ shows that they coincide
on $f^{-n}\xi^u | \xi^u(y)$.
Since $$\bigvee_{n=0}^\infty f^{-n}\xi^u=\epsilon$$ we
 conclude that $\nu(\cdot|\xi^u(y))=\eta|\xi^u(y)\equiv\mu^+(\cdot|\xi^u(y))$.
 
Since $f^n\xi^u$ is also the Pesin partition for any $n$, the conditional
measures of $\nu$ and $\mu^+$ coincide on it. Passing to the limit as
$n\to\infty$, we get the required agreement of $\nu$ and $\mu^+$. \qed

\give Proof of Theorem 3.1.
By the ergodic theorem $\nu$-almost every point $p$ is equidistributed with
respect to $\nu$, that is
 $$\lim_{n\to\infty}{1\over n}\sum_{i=0}^{n-1}\phi(f^n(p))=\int\phi \nu
\eqno (3.3)$$
holds for any continuous function $\phi$ on ${\bf C}^2$ with compact support.
Let $\nu_x=\nu(\cdot|\xi^u(x))$ denote the conditional measure on the
Pesin piece $\xi^u(x)$.
Then for almost every  $x$ we have that $\nu_x$-almost every point
in $\xi^u(x)$ is equidistributed with respect to $\nu$.
By bounded convergence we can average (3.3) over $\xi^u(x)$:
$$\eqalign{\int\phi \nu&=\lim_{n\to\infty}\int{1\over n}
\sum_{i=0}^{n-1}\phi(f^i(p))\nu_x(p)\cr
&=\lim_{n\to\infty}\int\phi(p) ({1\over n}\sum_{i=0}^{n-1}f_*^i(\nu_x))(p)\cr}$$
    Since this holds for any continuous function we have:
$${1\over n}\sum_{i=0}^{n-1}f_*^i(\nu_x)\to \nu \eqno (3.4)$$
in the weak topology of measures.
 
 On the other hand, let $\mu^+_x = \mu^+(\cdot|\xi^u(x))$ denote the
 normalized measure $\mu^+|\xi^u(x)$. Then it follows from [BS4] that
$$f_*^i(\mu^+_x)\to \mu. \eqno (3.5)$$
To see this we observe that the proof of Lemma 4.1 in [BS4] is still valid if
 $L$
is a continuous function and $S$ is a current for which the conclusion of
 Theorem
1.6 holds. We take $L=G^+$ and $S$ to be the current of
integration of the set $\xi^u(x)$. In the notation of [BS4] we have
$\Theta_n=\mu^+$ and $\nu_n=d^n\mu^+_x$. The conclusion of Theorem 1.6 holds for
$S$ by the third remark following the proof of Theorem 1.6 (cf. example
 following
Corollary 1.7).
 
Since $\nu_x = \mu^+_x$ by Proposition 3.2,  properties (3.4) and (3.5)
 yield $\nu=\mu$.  This completes the proof of the theorem. \qed
 
\proclaim Corollary 3.3. The conditional projective measure class
$\bar\mu(\cdot|W^u(x))$
is induced by the current $\mu^+$.
 
\give Remark.  The Jacobian of $f$ with respect to the conditional measures on 
 the
unstable manifolds is thus $d$.  This can be considered as the natural analogue
of the balanced property of the Brolin measure.
\medskip
 It is known that for a polynomial endomorphism $P$ of the complex plane
the characteristic exponent $\chi$ of the measure of maximal entropy (which
coincides with harmonic measure of the Julia set $J(P)$)
is greater or equal than $\log d$. Moreover, $\chi=\log d$ if and only if $J(P)$
is connected (see [Man], [Pr]). Here we discuss related
properties of polynomial automorphisms of $\cx2$.
 
It follows from the Lai Sang Young formula (see \S2.5)
 and Corollary 3.3 that for $\mu$~a.e. $x$,
 
$$HD(\mu^-|_{W^u(x)})={h_\mu(f)\over|\chi^s|}={\log d\over|\chi^s|}.\eqno(3.6)$$

Let us consider for a moment the dissipative case, i.e.\ $|a|<1$,
 where $a$ is
the (constant) Jacobian determinant of $f$.  We have
$|\chi^s|=\chi^u-\log|a|$.  It was shown in [BS4] that $\chi^u\ge\log d$,
 so in
this case, Young's formula gives $$HD(\mu^- |_{W^s(x)})<1\eqno(3.7)$$
for $\mu$~a.e.\ $x$.  By Corollary 3.3 and
the fact that the conditional measures are in the measure class of  harmonic
measure, we have:
 
\proclaim Corollary 3.4.  If $f$ is dissipative,
 then for $\mu$~a.e.~$x$ the harmonic
measure of $W^s(x)\cap J^-$ inside $W^s(x)$ has Hausdorff dimension strictly
less than 1.
 
In the following result, we relate the topological property of the
connectedness of $J$ to the rate of expansion of $f$. 
 
\proclaim Theorem 3.5.  If the map $f$ is hyperbolic, and if $J$ is connected,
then $\chi^u=\log d$ and $\chi^s=\log|a|-\log d$.
 
\give Proof.  Since $f$ is hyperbolic, $J$ has a local product structure at any
point $p$.  That is, there are neighborhoods $V^s$ of $p$ in $W^s(p)$ and $V^u$
in $W^u(p)$  such that $(V^s\cap J^-)\cap(V^u\cap J^+)$ is homeomorphic to a
neighborhood of $p$ in $J$.
 
We claim that either $W^s(p)\cap J^-$ or $W^u(p)\cap J^+$ has the property: 
There is a neighborhood $U$ of $p$ such that every connected component of
$W^s(p)\cap J^-\cap U$ (resp.\ $W^u(p)\cap J^+\cap U$) is noncompact.  For
otherwise there are compact connected components which are arbitrarily small and
arbitrarily close to $p$ inside both $W^s(p)\cap J^-$ and $W^u(p)\cap J^+$.   Thus the
product neighborhood of $p$ in $J$ contains arbitrarily small compact, connected
components. But in this case, $J$ is not connected, which proves the claim.
 
Thus we may assume that $W^u(p)\cap J^+$ has this property.  If $q\in J$ is close
to $p$, then by the local product structure, $W^u(q)\cap J^+$ also has this
property.  It follows that $U\cap W^u(q)-J^+$ is simply connected.  By a Theorem
of Makarov [Mak], the harmonic measure of $W^u(q)\cap J^+$ has Hausdorff
dimension 1. Now since this holds for a set of $q$ of positive measure, we
conclude from the formula of Young, with the stable manifolds replaced by the
unstable manifolds, that $\chi^u=\log d$.\qed

\bigskip

\section 4.  Product Structure of $\mu$

In this Section we will show that there are sets (``Pesin boxes")
 on which $\mu$ has a local
product strucure, and the union of these sets has full $\mu$ measure.  The main
step in doing this is to study the holonomy map along the stable/unstable
manifolds and to show that the conditional measures of $\mu$ are preserved by
the holonomy map.
 
We consider a family $\cM$ of complex manifolds.  A complex manifold $D$ is a
{\it transversal} to $\cM$ if $D$ intersects each $M\in\cM$ in a unique point,
and this intersection is transverse.  Let $D_1$ and $D_2$ be two transversals
to $\cM$, and set
$$X_i=\bigcup_{M\in\cM} D_i\cap M$$
for $i=1,2$.  We define the {\it holonomy map}
$$\chi:=\chi({D_1,D_2,\cM}):X_1\to X_2$$
as $\chi(x_1)=M(x_1)\cap D_2$, where $M(x_1)\in\cM$ is the unique
manifold containing $x_1$.
 
Throughout this Section we will consider the case where $\cM$ is a family of
stable (or unstable) disks given by the Pesin theory (see the discussion in
\S2).  If  $r>0$ is sufficiently small, and if $x\in\cR$ satisfies $r(x)\ge r$,
then each stable (or unstable) disk $W^s_r(x)$ (or $W^u_r(x)$) is given as a
graph over the  $r$-ball in the tangent space $E^{s}(x)$ (or $E^u(x)$).  More
generally, we will work with complex disks that are graphs over $E^{u}_r(x)$,
i.e. which have the form  $$M=\{(z,\varphi(z)): z\in
E^u_r(x),\varphi(z)\in E^s_r(x)\}.$$
We will say that two such graphs are {\it $C^1$ close} if their corresponding
graphing functions $\varphi$ are $C^1$ close.

For any subset $F\subset \{x\in\cR:r(x)\ge r\}$, we write
$$W^{s/u}_{loc}(F)\equiv W^{s/u}_r (F)=\bigcup_{x\in F} W^{s/u}_r(x)
 {\rm\ \ and\ \ }
\cF^{s/u}=\{W^{s/u}_r(x):x\in F\}.\eqno(4.1)$$ For $0<\kappa<1$ and for
$x_0\in\cR$ with $r(x_0)\ge r$, we let $F=\{x\in\cR\cap B(x_0,\kappa
r):r(x)\ge r\}$.  We may choose $\kappa$ sufficiently small that if $D_1$ and
$D_2$ be disks such that $D_j$ is within $\kappa r$ (in the $C^1$ topology) of
$W^u_r(x_j)$ for some $x_j\in F$, then $D_j$ is transversal to the family
$\cF^s$ for $j=1,2$.  It follows that the holonomy map $\chi(D_1,D_2,\cF^s)$ is
defined.

\def\rs{{r_0}}
For $m_0>0$, $\rs=r/8>0$, and $x\in\cR$ such that $r(x)\ge r$, we consider the
property $$\mu^+|_{W^u_r(x)}(W^u_\rs(x))\ge m_0.\eqno(4.2)$$
For $C<\infty$ we also consider the properties
$${\rm dist}(f^n(x),f^n(y))\le
Ce^{-n\lambda} {\rm\ \ for\ }n\ge1{\rm\ and\ }y\in W^s_r(x)\eqno(4.3)$$
$${\rm dist}(f^{-n}(x),f^{-n}(y))\le
Ce^{-n\lambda} {\rm\ \ for\ }n\ge1{\rm\ and\ }y\in W^u_r(x)\eqno(4.4)$$
Choosing $m_0>0$ sufficiently small and $C<\infty$ sufficiently large, we
have
$$\mu(J-Q)<\epsilon\eqno(4.5)$$
where
$$Q=\{x\in\cR:r(x)\ge r,{\rm\ and\ }(4.2), (4.3),(4.4){\rm\
hold}\}.$$
 
Now let
$$S=\{x\in J:f^n(x)\in Q{\rm\ for\ infinitely\ many\ }n\}.\eqno(4.6)$$
By the Ergodic Theorem, $\mu(S)=1$.  In the sequel, we let $Q$ denote the set
$Q\cap S$, which differs from the original $Q$ by a set of measure zero.  Now
 let
us fix $x_0\in Q$ and use the following notation:
$$F=Q\cap B(x_0,\kappa r),\eqno(4.7)$$
let $D=W^u_r(x_0)$, and let $D^\prime$ be a transversal which is
within $C^1$-distance $\kappa r$ of $D$.  The domain of the holonomy
map $\chi(D,D^\prime,\cF^s)$ is given by $$X=D\cap W^s_r(F)$$
and the range is
$$X^\prime=D^\prime\cap W^s_r(F).$$
 
We recall that the construction of the Pesin unstable manifolds (as given,
for instance, in [PS]) may be carried out by applying the graph transform,
starting with disks, called ``trial disks," that are transverse to the stable
direction.  It is shown that these trial disks, under forward iteration, approach
the stable manifolds in a semi-global $C^1$ sense.  Now let us consider a large
$n$ such that $f^nx_0\in Q$.  We define $y_0=D^\prime\cap W^s_r(x_0)$ and view
$D^\prime$ as a trial disk for the unstable manifold $D=W^u_r(x_0)$.  Let
$D_n^\prime$ denote the portion of $f^nD^\prime$ which can be represented as a
graph over $E^u_r(f^nx_0)$ and which contains $f^ny_0$.  By [PS, Corollary 3.11]
$D_n^\prime$ converges to $D_n\equiv W^u_r(f^nx_0)$ in the $C^1$ topology, and
the  distance is bounded by $Ce^{-n\lambda}$.  Let us define $$X_n=(f^nX)\cap
W^u_{\kappa r}(f^nx_0).$$ Evidently, $f^{-n}X_n\subset X$.
 
For each $x\in f^nX_n$, the local stable manifold $W^s_{r(f^nx)}(f^nx)$
intersects $D_n^\prime$ transversally because
$r(f^nx)\ge(1+\epsilon)^{-n}r$, and the
angle between $W^s_{r(f^nx)}(f^nx)$ and $D_n$ is at least
$\theta_0(1+\epsilon)^{-n}$, whereas $D_n^\prime$ is exponentially close to
$D_n$.  Thus $D_n^\prime$ and $D_n$ are transversals to the family
$\cF_n=\{W^s_{r(f^nx)}(f^nx):x\in X_n\}$, and so the holonomy
$$\chi_n:=\chi(D_n, D_n^\prime, \cF_n)\eqno(4.8)$$ is defined.
 
\proclaim Lemma 4.1.  If $n$ is sufficiently large, and $f^nx\in Q$, then
$\chi_n(f^nx)=f^n\chi(x)$
for all $x\in X_n$.  Given $\rs<\kappa r$, $n$ may be taken sufficiently large
that for $a\in f^{-n}X_n$
$$\eqalign{
\chi_n (X_n\cap B(f^na,\rs-2Ce^{-n\lambda}))\subset
(\chi_n X_n)&\cap B(f^n\chi (a),\rs)\subset\cr
\subset&
\chi_n (X_n\cap B(f^na,\rs+2Ce^{-n\lambda})).\cr}$$
 
\give Proof.  For $x\in X_n$, let $y=\chi(x)$, and let $\gamma$ be a path
inside $W^s_r(x)$ connecting $x$ to $y$.  Then $f^n\gamma$ lies inside
$f^nW^s_r(x)$.  Further, by (4.3), $f^n\gamma$ has diameter less
than $Ce^{-n\lambda}$, and thus $f^n\gamma\subset W^s_{r(f^nx)}(f^nx)$.  Since
$f^n\gamma$ connects $f^nx$ to $f^ny$ inside $W^s_r(f^nx)$, it follows that
$\chi_n(f^nx)=f^ny$.
 This proves the first assertion.  The required inclusions are now 
a consequence of (4.3).\qed
 
\give Remark. Sometimes abusing rigour we will write
$\chi\circ f^n=f^n\circ \chi$ and say that $f^n$ commutes with holonomy. Let us
use the notation  $y_0=\chi(x_0)$, $x_n=f^n(x_0)$, and $y_n=f^n(y_0)$. 
 
\proclaim Lemma 4.2.  Let $\{\eta_n\}$ be a sequence of numbers decreasing to
zero.  Let us pass to a subseqence $n=n_j$ for which $x_n\in Q$, and let
$D_n^\prime$
be a sequence of complex disks such that $dist_{C^1}(D_n,D_n^\prime)\le \eta_n$.
Then there exists $\rho$ with $\rs/2\le\rho\le\rs$ such
that
$${\lim_{n\to\infty}}^\prime\left[\mu^+|_{D_n^\prime}B(y_n,\rho\pm2Ce^{-n\lambda
 })
- \mu^+|_{D_n}B(x_n,\rho\mp 2Ce^{-n\lambda})\right]=0,$$
where ${\lim}^\prime$ means that the limit is taken through a further
subsequence.
 
\give Proof.  Without loss of generality, we may assume that $Q$ is compact,
and a subsequence of $\{x_n\}$ converges to $\bar x\in Q$.  Thus the unstable
disks $D_n$ converge in $C^1$ to $\bar D=W^u_r(\bar x)$.  Now choose
$\rs/2\le\rho\le \rs$ such that $\mu^+|_{\bar D}$ puts no mass on
$\partial B(\bar x,\rho)$.  The Lemma then follows because the measures
$\mu^+|_{D_n^\prime}$ converge weakly to $\mu^+|_{\bar D}$.\qed
 
\proclaim Lemma 4.3.  If, in addition to the hypotheses of Lemma 4.2, we
require that $x_n\in Q$, then
$${\lim_{n\to\infty}}^\prime\mu^+|_{D_n^\prime}B(y_n,\rho\pm2Ce^{-n\lambda})
\left(\mu^+|_{D_n}B(x_n,\rho\mp 2Ce^{-n\lambda})\right)^{-1}=1.$$
 
\give Proof.  Lemma 4.3 follows from Lemma 4.2 by property (4.2).\qed
 
\proclaim Lemma 4.4.  Let $F\subset Q$, $D=W^u_r(x_0)$, and
$D^\prime$ be as above.  With the notation $\nu:=\mu^+|_{D}$,
$\nu^\prime=\mu^+|_{D^\prime}$, and $\chi=\chi(D,D^\prime,\cF^s)$, we have
$$\chi_*(\nu |_{X })=\nu^\prime|_{X^\prime}.\eqno(4.9)$$
 
\give Proof.  It will suffice to show that (4.9) holds for $X $ replaced by
$X \cap B(x,\epsilon)$ for some small $\epsilon>0$.  Then we can add over a
partition of $X $ to obtain (4.9).  We will define two coverings $\cC^\pm$ of
$X $ and a covering $\cC^\prime$ of $X^\prime$.  The coverings will have the
property that if $a\in X $, there are elements $C^\pm(a)\in\cC^\pm$ and
$C^\prime(\chi(a))\in\cC^\prime$ of arbitrarily small size containing $a$ such
that  $$\chi C^-(a)\subset C^\prime(\chi(a))\subset\chi C^+(a).\eqno(4.10)$$
For $a\in X $ we may choose $n$ arbitrarily large such that $f^na\in Q$.
  We define
$$C^\pm_n(a)=f^{-n}(W^u_r(f^na)\cap B(f^na,\rho\pm2Ce^{-n\lambda}))$$
with $\rho$ as in Lemmas 4.2 and 4.3.  In analogy with notation used earlier in
this Section, we let $D_n^\prime$ denote the portion of $f^nD^\prime$
 which lies as a graph
over $E^u_r(f^n a^\prime),\; a^\prime=\chi(a),$
 and which contains $f^n\ a^\prime$.  Now we define
$$C^\prime_n( a^\prime)=f^{-n}(D_n\cap B(f^n a^\prime,\rho)).$$
The inclusions in (4.10) are a consequence of Lemma 4.1.
By Lemma 4.3, we have
$${\lim_{n\to\infty}}^\prime {\nu (C^\pm_n(a))\over \nu^\prime(C_n^\prime
(a^\prime))}=1.\eqno(4.11)$$
 
By the overflowing property of the unstable disks, $f^{-n}:W^u_r(f^na)\to
W^u_r(a)$.  Since $W^u_r$ is a graph, we may identify it with the disk
$\{|\zeta|<r\}$, and thus we may consider $h(\zeta):=f^{-n}(r\zeta)$ as a
univalent mapping of the disk $\{|\zeta|<1\}$ to $\cx{}$.  By the Koebe
Distortion Theorem,
$$\left|{f^{\prime\prime}(\zeta)\over f^\prime(\zeta)}\right|\le
{(1+\rs/r)^4\over(1-\rs/r)^5}$$
for $|\zeta|<\rs/r$.
The image of the disk $\{|\zeta|<\rho\}$ is a convex set if
$(1+\rs/r)^4(1-\rs/r)^{-5}\le(2\rho)^{-1}$.  We conclude, then, since
$C^\pm_n(a)$ is the image of such a disk, and since $\rho\le\rs=r/8$, that
$C^\pm_n(a)$ is convex.  Similarly, $C^\prime( a^\prime)$ is convex.
 
Now let $E \subset X $ be a compact subset, and let $E^\prime=\chi E $.  For
$\delta>0$, choose an open set $\cO\subset W^u_r(x)$ containing $E $ such
that
$$\nu (\cO)<\nu (E )+\delta.$$
 
The coverings $\cC^\pm$ and $\cC^\prime$ are fine in the sense that any point
is contained in an element of arbitrarily small diameter.  Since the elements
of the cover are convex, we may apply the Covering Theorem of A.P.~Morse [Mo]
to conclude that there is a disjointed family $\{C^\prime_j:
j=1,2,\dots\}\subset\cC^\prime$ such that
$$\nu^\prime\left(E^\prime-\bigcup_{j=1}^\infty
C^\prime_j\right)=0.\eqno(4.12)$$ Each $C^\prime_j$ is of the form
$C^\prime_{n_j}(\chi(a_j))$ for some $a_j$ and $n_j$.  The corresponding sets
$C^-_j$ in the cover $\cC^-$ satisfy $C^-_j\subset\chi^{-1}(C^\prime_j)$ by
(4.10) and are thus pairwise disjoint.  Since $\chi^{-1}$ is continuous, and
since the diameters of the $C^\prime_j$ may be taken arbitrarily small,  we may
assume that $C^-_j\subset\cO$.  Thus $$\nu (\cO)\ge\sum_{j=1}^\infty\nu
(C^-_j).$$ Since we may take the diameters arbitrarily small, it follows from
(4.11) that
$$\nu(\cO)\ge(1+\delta)^{-1}\sum_{j=1}^\infty\nu^\prime(C^\prime_j).$$ By
(4.12), then, $\nu(\cO)\ge\nu^\prime(E^\prime)$.  It follows that $$\nu (E
)\ge\nu^\prime(E^\prime).$$
 
Now if we cover $E $ by a disjointed subcover of $\cC^+$ and repeat the
previous argument, we conclude that
$$\nu (E )\le\nu^\prime(E^\prime).$$
Thus $\nu(E)=\nu^\prime(E^\prime)$, and this completes the proof.\qed

\proclaim Theorem 4.5.  Let $F\subset Q$, and $\cF^s$ be as
above.  Let $D_1$, $D_2$ be two transversals, and set
$\mu_j:=\mu^+|_{D_j}$ for $j=1,2$,  $X_j=D_j\cap W^s_r(F)$.  Then the
holonomy $\chi:=\chi(D_1,D_2,\cF^s)$ satisfies
$$\chi_*(\mu_1|_{X_1})=\mu_2|_{X_2}.$$
 
\give Proof.  We may assume that for each $x_1\in X_1$, there is a point
$z\in Q$ such that $x_1,\chi(z)\in W^s_{\kappa r}(x^\prime)$.
 For otherwise we may apply $f^n$ and use (4.3) and the fact
that $f^n\circ\chi=\chi\circ f^n$.
 
As in Lemma 4.4 we work locally on $X_1$, so we may assume that
$W^u_r(z)$ is a transversal to $\cF^s$.  Let us define
$\chi_1=\chi(D_1,W^u_r(z),\cF^s)$ and
$\chi_2=\chi(W^u_r(z),D_2,\cF^s)$.  By Lemma 4.4, then,
$\chi_2\circ\chi_1=\chi$ takes $\mu_1|_{X_1}$ to $\mu_2|_{X_2}$.\qed
 
Let $F\subset Q$ denote a compact subset, and let $W^s_r (F)$,
${\cal F}^{s/u}$ be as in (4.2).
 If the diameter $\delta$ of $F$ is sufficiently
small, then we may assume that $F$ is contained in a $\delta$-ball about the
origin, and that every leaf $W^s_{r}(x)$ (resp. $W^s_{r}(x)$) is a graph
over the horizontal (resp. vertical) coordinate axis.  Further, for $x\in F$,
$W^s_{r}(x)$ is transversal to ${\cal F}^u$, and $W^u_{r}(x)$ is
transversal to ${\cal F}^s$.  Since the holonomy induces a homeomorphism on
transversals, there is a fixed compact set $P^u$ which is homeomorphic to
$W^u_{r}(x)\cap W^s_r(F)$ for all $x\in F$.  Similarly, there is a fixed
$P^s$ which is homeomorphic to
$W^s_{r}(x)\cap W^u_r(F)$ for all $x\in F$.  We call the set
$P:=W^s_r(F)\cap W^u_r(F)$ the {\it Pesin box} generated by $F$, and we note
that $ P$ is naturally homeomorphic to $P^s\times P^u$.  It is evident
that, up to a set of measure zero, ${\cal R}$ is a countable union of (not
necessarily disjoint) Pesin boxes.
 
If $ P$ is a Pesin box, then the partitions ${\xi}^{s/u}$ of $P$,
whose elements are $W^{s/u}_{r}(x)\cap P$ are measurable partitions of
$P$.  We let $c:=\mu( P)$ so that $\nu:=c^{-1}\mu\contract\tilde
P$ is a probability measure.  It follows from Theorem 3.1 that the conditional
measures of $\nu$ are given by
$$\nu(\cdot|\xi^{s}(x))=c^{s}(x)^{-1}\mu^{-}|_{W^{s}_{r}(x)}\contract P$$
where $c^{s}(x)$ is the total mass of $\mu^{-}|_{W^{s}_{r}(x)}\contract P$,
 and a similar expression for $\nu(\cdot|\xi^{u}(x))$.
By Theorem 4.5 we see that $c^s=c^{s}(x)$ is constant for $x\in F$.  In fact:

\proclaim Theorem 4.6.  If $ P$ is a Pesin box as above, then the
holonomy maps along ${\cal F}^{s/u}$ preserve the conditional measures of
$\nu=c^{-1}\mu\contract P$.
 
To explore the product structure further, we let $E^{s/u}\subset P^{s/u}$ be
Borel sets.  For $x\in F$, we may define the measures $\lambda^{s/u}$ on
$P^{s/u}$ to be the measures induced by  the conditional measures
$c^{s/u}\nu(\cdot|\xi^{s/u}(x))$ via the homeomorphism between $\xi^{s/u}$
and $P^{s/u}$.   By Theorem 4.6, the measures $\lambda^{s/u}$ are independent of
the point $x\in F$.  By properties (ii) and (iii) of conditional measures, we
have   $$\eqalign{
\nu(E^s\times P^u)&=\int_{x\in\tilde
P}\nu(E^s|\xi^s(x))\nu(x)\cr
&=\int_{x\in\tilde
P}(c^s)^{-1}\lambda^s(E^s)=(c^s)^{-1}\lambda^s(E^s).\cr}$$
Similarly, we have
$$\eqalign{
\nu(E^s\times E^u)&=\int_{x\in P^s\times E^u}\nu(E^s|\xi^s(x))\nu(x)\cr
&=\int_{x\in P^s\times E^u}(c^s)^{-1}\lambda^s(E^s)\cr
&=(c^s)^{-1}\lambda^s(E^s)\nu(P^s\times
E^u)=(c^sc^u)^{-1}\lambda^s(E^s)\lambda^u(E^u).\cr}$$  Thus we have the
following.

\proclaim Theorem 4.7.  If $P$ is a Pesin box, then there are measures
$\lambda^{s/u}$ on $P^{s/u}$ such that $\mu\contract
P=\lambda^s\otimes\lambda^u$ has the structure of a product measure.
 
\proclaim Corollary 4.8. The measure $\mu$ is Bernoulli.
 
\give Proof.  D.\ Ornstein and B.\ Weiss discuss invariant
measures with nonzero Lyapunov exponents in [OW, p.\ 86].  Given such a measure
which is mixing with respect to $f$, they remark that it is Bernoulli if it
is locally equivalent to a product measure with respect to the stable and
unstable manifolds.  By Theorem 4.7, then, we conclude that $\mu$ is
Bernoulli. \qed
 
Since the entropy, $\log d$, depends only on the degree of $f$, and the
entropy is the unique invariant for Bernoulli measures, it follows that any
two polynomial automorphisms with the same degree are measurably conjugate
with respect to their equilibrium measures.
 
\bigskip
\bigbreak

\section 5. Uniformly Laminar Currents

Let  $\Omega\subset\cx n$ be an open set, and let $\cD^{p,q}$ denote the smooth
 $(p,q)$-forms
$\alpha=\sum\alpha_{IJ}dz^I\wedge d\bar z^J$, $|I|=p$, $|J|=q$, with compact
 support in
$\Omega$.  The dual space $\cD_{p,q}$ of $\cD^{p,q}$ is the set of $(p,q)$-{\it
 currents}
or {\it currents of bidimension $(p,q)$}.  A current of dimension 0 acts on test
functions and may thus be considered as a distribution.  $\cx n $ itself may be
 identified with
the $2n$-dimensional current $[\cx n]$, which acts on an $(n,n)$ form $\varphi$
 by
integration: $[\cx n](\varphi)=\int\varphi$.  If $T$ is a $(p_1,q_1)$-current,
 and $\psi$ is a
smooth $(p_2,q_2)$-form, then the {\it contraction} $T\contract\psi$, defined by
$$(T\contract\psi)(\varphi)=T(\psi\wedge\varphi)$$ is a
 $(p_1-p_2,q_1-q_2)$-current.  The space
$\cA^{n-p,n-q}$ of smooth $(n-p,n-q)$ forms on $\cx n$ may be identified with a
 set of currents
of bidimension $(p,q)$ via the mapping
$$\cA^{n-p,n-q}\ni\psi\mapsto[\cx n]\contract\psi\in\cD_{p,q}.$$
 
The {\it mass norm} of a current $T$ is given by
$${\bf M}[T]=\sup_{|\varphi|\le1}|T(\varphi)|.$$
If $T$ is an $(0,0)$ current, then the mass norm is finite if and only if $T$ is
 represented
as a distribution by a finite, signed Borel measure $\nu$, and ${\bf M}[T]$ is
 the total
variation of $\nu$.   A current $T$ is {\it representable by integration} if
 $\chi T$ has
finite mass norm for any test function $\chi$ on $\Omega$.  If $T$ is
 representable by
integration, then there is a Borel measurable function $t$ from $\Omega$ to the
$(p,q)$-vectors (the dual of the $(p,q)$-forms) and a Borel measure $\nu$ on
 $\Omega$
such that $T=t\nu$ holds in the sense that
$$T(\varphi)=\int_{x\in\Omega}\langle\varphi(x),t(x)\rangle\nu(x).$$
We will require that $|t|^*=1$ at $\nu$ a.e.~point.  ($|\cdot|^*$ denotes the
 norm
on $(p,q)$-vectors which is dual to the norm on $(p,q)$ forms.)  In this case
 $t$
and $\nu$ are uniquely determined, and $\nu=|T|$ is the variation measure
associated with the current $T$.  We will call $t\nu$ the {\it polar
representation} of $T$.  If $T$ is representable by integration, and if
$S\subset\Omega$ is a Borel subset, then we will use the notation  $$T\contract
S=t\nu\contract S$$ for contraction, which coincides with restriction in this
case.

A $(p,p)$-current $T$ is {\it positive} if
$T(i\alpha_1\wedge\bar\alpha_1\wedge\ldots\wedge
 i\alpha_p\wedge\bar\alpha_p)\ge0$ for all
$(1,0)$ forms $\alpha_j=\sum_k\alpha^k_jdz_k$ with compact support.  This
 definition of
positivity is analogous to the positivity of a distribution.  And as in the case
 of
distributions, a positive current is representable by integration.  Further, if
 we let
$\beta=\sum{i\over 2}dz_j\wedge d\bar z_j$ denote the standard K\"ahler form on
 $\cx n$, then
for a positive $(p,p)$ current $T$, the contraction $T\contract\beta^{p}/p!$ is
 a Borel
measure, and the mass norm ${\bf M}[T]$ is just the total variation of this
 measure.
 
Let $M$ be a $k$-dimensional complex manifold of $\Omega$.  If either $M$ is
 locally closed
(without boundary) or if $M$ is a smooth submanifold-with-boundary (or more
generally, if the area of $M$ is locally finite), then the pairing with test
$(k,k)$-forms given by  $$[M](\varphi)=\int_M\varphi$$
defines $[M]$ as a current of bidimension $(k,k)$ on $\Omega$.  We call $[M]$
 the {\it
current of integration} associated to $M$.  The mass norm of $[M]$ is the
 Euclidean
$2k$-dimensional area of $M$.  It is evident that $[M]$ is representable by
 integration, and
$$[M]=t_M\sigma_M,\eqno(5.1)$$
where $t_M$ is the $2k$-vector of norm 1 defining the
tangent space to $M$ (a vector which is uniquely defined, since $M$ is an
 oriented
submanifold of $\cx n$), and $\sigma_M=\cH^{2k}\contract M$ is the Hausdorff
$2k$-dimensional measure restricted to $M$.   The {\it boundary} $\partial T$ of
 a current
$T$ is defined by $$\partial T(\varphi)=T(d\varphi).$$
If $\partial M$ is regular, we may apply Stokes' theorem to obtain
$\partial[M](\xi)=\int_{\partial M}\xi$.  We say that
$T$ is {\it closed} if $\partial T=0$, and so $[M]$ is closed if $M$ has no
 boundary.
 
More generally, if $V$ is a (closed) subvariety of $\Omega$, then the set
 $Reg(V)$ of
regular points (where $V$ is locally a manifold) are a dense open set, and it
 may be shown
that $[V](\varphi)=\int_{Reg(V)}\varphi$ defines a positive, closed current.
 The device of
studying the current of integration $[V]$ has been useful in the study of metric
 properties
of $V$, such as the area growth.  For instance, the fact that $[V]$ is a current
 at all
corresponds to the fact that the area of $[Reg(V)]$ is locally bounded near
 singular
points.  And $\partial[V]=0$ holds because the amount of mass in a neighborhood
 of the
singular set is small.
 
It is useful to apply similar considerations to the stable and unstable
 manifolds.  However,
since $W^s(x)$ (resp. $W^u(x)$) is often dense in $\cW^s$ (resp. $\cW^u$) an
 individual
stable manifold does not define a current of integration, since the amount of
 mass is not
locally bounded.  Thus we wish to consider the whole stable and unstable
 laminations as
currents, as was suggested by Ruelle and Sullivan [RS] and Sullivan [S].
 
Let us consider a family of graphs of analytic functions $f_a:\Delta\to\Delta$,
$a\in\Delta$.  We assume that the graphs $\Gamma_a=\{(x,f_a(x)):x\in\Delta\}$
 are pairwise
disjoint, i.e. if $a_1\ne a_2$, then $f_{a_1}(x)\ne f_{a_2}(x)$ for all
 $x\in\Delta$.  We
denote the set of graphs as $\cG=\{\Gamma_a:a\in A\}$.  Without loss of
 generality, we may
take the parameter space to be a closed subset of the unit disk, and we may take
 $a=f(0)$.
Further, since the graphs are disjoint, it follows that $a\mapsto f_a$ is
 continuous.
 
A current $T$ on $\Delta^2$ is {\it uniformly laminar} if it has the form
$$T=\int_{a\in A}\lambda(a)\,[\Gamma_a]\eqno(5.2)$$
where $\lambda$ is a positive measure on $A$, the parameter space for the set
 $\cG$ of
graphs.  The action on a (1,1) form $\varphi$ is given by
$$T(\varphi)=\int_A\lambda(a)\int_{\Gamma_a}\varphi.$$
We say that a current $S$ is {\it
locally uniformly laminar} on an open set $\Omega$ if for each $p\in\Omega$
 there is a coordinate
neighborhood equivalent to $\Delta^2$ on which $S$ is uniformly laminar.  The
 currents of integration
$[\Gamma_a]$ are positive, closed currents on $\Delta^2$, so $T$, too, is
 positive and
closed.
 
For a transversal $M$ to the family $\cG$, the set
of all intersection points, $A_M$, could equally well be taken as a parameter
 space.
Further, let $M_1$ and $M_2$ be transversals.  Then the holonomy
map $\chi_{M_1,M_2}:A_{M_1}\to A_{M_2}$ gives a homeomorphism between parameter
 spaces.  For a
point $p\in\cx2$, we let $[p]$ denote the 0-current which puts a unit mass at
 the point $p$.
For each transversal, the current (measure) $[\Gamma_a\cap M]$ depends
 continuously on $a$.
We define the {\it restriction} of $T$ to $M$ by
 $$T|_M=\int_A\lambda(a)[\Gamma_a\cap
M],\eqno(5.3)$$ which is a measure on $M$.  If $M_1$ and $M_2$ are transversals,
 then the
restrictions are preserved by the holonomy map $\chi=\chi_{M_1,M_2}$, i.e.
$$\chi_*T|_{M_1}=T|_{M_2}.\eqno(5.4)$$
A family of measures $\{T|_M\}$ on transversals induces a {\it transversal
measure} on $\cW^s$ if it satisfies (5.4).
$T$ may be reconstructed from any transversal (or, equivalently, from any
family of transversal measures) as $$T=\int_{a^\prime\in
A_M}T|_M(a^\prime)[\Gamma_{a^\prime}].\eqno(5.5)$$ Equations (5.4) and (5.5) are
trivial if $T=[\Gamma_a]$ is a current of integration, and the general case is
obtained by integrating with respect to $\lambda$.  Let $h$ be a holomorphic
function on $\Delta^2$ such that $M=\{h=0\}$ and $dh\ne0$ on $M$.  Then
$\log|h|$ is locally integrable on each $\Gamma_a$, and
$$dd^c\log|h|[\Gamma_a]=[M\cap\Gamma_a]$$ holds in the sense of currents.  Thus
$${1\over 2\pi}T|_M=dd^c(\log|h|T).$$
 
We may ask, more generally, which positive, closed currents on $\cx2$ may be
 represented in
the form
$$T=\int_{a\in A}\eta(a)\,[V_a]\eqno(5.6)$$
where $A\ni a\mapsto V_a$ is a measurable family of varieties in $\cx2$, and
 $\eta$ is a
Borel measure on $A$.  This is closely related to the Choquet representation of
 $T$ as an
integral over extremal rays on the cone of positive, closed currents.  It is
 known that an
irreducible subvariety $V_a\subset\cx2$ generates an extreme ray (see [D] and
[L]).
 On the other
hand, not all extremal rays are of the form $c[V]$.  This
will also be a consequence
 of the
examples below.
 
\give Examples.  Let $(x,y)$ denote coordinates on $\cx2$, and define
$$u_1=\log^+|(x,y)|=\max\{0,{1\over2}\log(|x|^2+|y|^2)\}$$
$$u_2=\max\{\log|x|,\log|y|,0\}.$$
For $\alpha\in\cx2$, we let $L_\alpha$ denote the complex line through 0 and
 $\alpha$, and
we set $L_\alpha^+=L_\alpha\cap(\cx2-\bar{\bf B}^2)$.  Then we may compute
$$T_1:=dd^cu_1=2\pi\int_{\alpha\in{\bf P}^1}[L^+_\alpha]\sigma(\alpha) + S_1,$$
where $\sigma$ is normalized spherical measure on ${\bf P}^1$, and $S_1$ is
 supported on
$\partial{\bf B}^2$.  Similarly,
$$\eqalign{T_2:=dd^cu_2&=\int_0^{2\pi}[L^+_{(1,e^{i\theta})}]\,d\theta +
\int_0^{2\pi}[x=e^{i\theta},|y|<1]\,d\theta\cr
&+\int_0^{2\pi}[y=e^{i\theta},|x|<1]\,d\theta+S_2\cr}$$
where $S_2$ is supported on the 2-torus $\{|x|=|y|=1\}$.
 
It is evident, then, that $T_1$ is locally uniformly laminar on
 $\cx2-\partial{\bf
B}^2$, and $T_2$ is locally uniformly laminar on $\cx2-\{|x|=|y|=1\}$.
 
Now if $T=t\lambda$ is any positive current satisfying $T\le T_1$, then at
 $\lambda$ a.e.
point $\alpha\in\cx2-\bar{\bf B}^2$, the (1,1) vector $t(\alpha)$ must be
 tangent to
$L_\alpha$.  If, in addition, $T$ has the form (5.6), then it follows that for
 $\eta$ a.e.
$a$ the  variety $V_a$ must be contained in $L_\alpha$ for some $\alpha$.  Since
 $V_a$ is a
subvariety, we must have $V_a=L_\alpha$.  On the other hand, since $T_1=0$ on
 ${\bf B}^2$,
it follows that $T=0$ on $\cx2-\partial{\bf B}^2$.  But now for $\eta$ a.e. $a$,
 we must
have $V_a\subset\partial{\bf B}^2$, which is impossible, so $T=0$.
A similar argument shows that if $0\le T\le T_2$ and $T$ has the form (5.6) then
 $T=0$.
 
These examples then show that: {\sl There are extreme rays in the cone of
 positive, closed
currents which are not generated by currents of integration over varieties.}
This observation was made by Demailly in [D], using the current $T_2$ written in
 a somewhat
different form.
 
Sullivan conjectured in [S] that a positive, closed current might be
written locally in the form (5.6) on a
dense, open set.  This cannot be the case, however, because of the following
 examples, which
are taken from [BT2].  For a number $r>0$ let $\chi_r(z)=rz$ denote
dilation, and for a point $a\in\cx2$ let $\tau_a(z)=z+a$ denote translation.
 Let
$\{r_j:j=1,2,3,\dots\}$ be dense in ${\bf R}^+$, and let $\{a_j:j=1,2,3,\dots\}$
 be dense in
$\cx2$.  Then the currents $$\tilde T_1=\sum 2^{-j}\chi_{r_j*}T_1\eqno(5.7)$$
$$\tilde
 T_2=\sum
2^{-j}\tau_{a_j*}T_2\eqno(5.8)$$ are positive and closed, and both have the
property of
 being nowhere
locally uniformly laminar.  From this it may be shown that neither current can
be represented in the form (5.6) on any open set.  

We note that the manifolds
of $\tilde T_1$ intersect correctly in the sense of \S6, although $\tilde T_1$
is not a weakly laminar current, even locally (cf.\ Proposition 6.2).  In fact,
if $L_j$ is uniformly laminar on an open set $U_j$, and if $\sum_{j=1}^\infty
 L_j\le
T_1$, then in fact $\sum_{j=1}^\infty L_j\le T_1-S_1$ (with a similar
property for $T_2$).  Thus  $${\bf M}[\sum L_j]\le {\bf M}[T_1]-{\bf M}[S_1],$$
and so $T_1$ and $T_2$
cannot be approximated from below by uniformly laminar currents, even in the
 sense of measure.

\give Remark.
Let us observe that it is possible to define the  wedge products $\tilde
 T_j\wedge\tilde T_j$ for
$j=1,2$ (see \S8.)  We do not know of an example as above with the additional
 property that
$\tilde T_j\wedge\tilde T_j=0$.
 
\bigbreak
\section 6.  Laminar Currents

The currents that arise in dynamical systems often derive their structure
from the stable and unstable manifolds.  The examples in \S5 show
that the category of positive, closed currents is too general for the dynamical
context.  Stable (or unstable) manifolds have no self-intersections and are
pirwise disjoint, so a represtentation (5.6) should involve the additional
requirement that the varieties $V$ be pairwise disjoint.  In fact, the context
in which currents have been constructed from dynamical systems has been the
uniformly hyperbolic case, and the currents obtained in this case are
uniformly hyperbolic.  In the case of a hyperbolic measure, this uniformity is
lost, and so we turn to the study of laminar currents.  The philosophy behind
the Sullivan conjecture is substantiated by Proposition 6.2 below, which
says: {\sl A laminar current is uniformly laminar outside a set of small
measure}.  
 
We say that two manifolds $M_1$ and $M_2$ {\it intersect
correctly} if either $M_1\cap M_2=\emptyset$ or $M_1\cap  M_2$ is an open subset
of $M_j$ for $j=1,2$, i.e.~they intersect in a set of codimension 0. We consider
a measurable set $A\subset\cx{}$ and a measurable function $f:\Delta\times
A\to\cx 2$ such that $f(\zeta,a)$ is an analytic injection in $\zeta$ for fixed
$a$.  We assume that any pair of image disks
$$M_a=\{f(\zeta,a):\zeta\in\Delta\}$$ intersects correctly.  Let $\lambda$
denote a $\sigma$-finite measure on $A$. If $$\int_A\lambda(a)\,{\bf M}[M_a\cap
U]<\infty\eqno(6.1)$$ for all relatively compact open sets $U\subset\cx2$, then
$$T=\int_{a\in A}\lambda(a)\,[M_a]$$ defines a positive current on $\cx2$.  A
current obtained in this way is called a {\it weakly laminar current} on $\cx
2$.  The current $T$ is {\it laminar} if the disks $M_a$ are pairwise disjoint. 
With suitable modifications, we can also define (weakly) laminar currents on an
open set $\Omega\subset\cx2$.  Thus if $U$ is open, then $T\contract U$ is again
(weakly) laminar.  We will say that $T$ is {\it represented} by the data
$(A,\cM,\lambda)$.  We note that for fixed $a\in A$, the function $\zeta\mapsto
f(\zeta,a)$ in the definition of $T$ is far from unique.  If we fix $[M_a]$,
then we can replace $f(\cdot,a)$ by any holomorphic imbedding
$f^\prime:\Delta\to M_a$  such that $M_a-f^\prime(\Delta)$ has zero area.
 
The parametrizing function $f$ in the definition  is not, strictly speaking,
necessary.  If we consider $\tilde M:=\bigcup_{M\in\cM} M $ be a
total space, then $\cM$ is a partition of $\tilde M$, and $A=\tilde M/\cM$ is
 the
quotient.  The essential point is the requirement that this partition be
measurable. We say that the families $\cM_1$ and $\cM_2$ {\it intersect
correctly} if all of the component manifolds intersect correctly.
 
A Borel set $E$ is a {\it
carrier} for $T$ if $T\contract E=T$, or equivalently,  $E$ carries all the mass
of $|T|$.  A carrier for a (weakly) laminar current may be taken to be a union
of complex disks.
 
\proclaim Lemma 6.1.  Let $T_j$, $j=1,2,3,\dots$ be a sequence of weakly laminar
currents with representations $(A_j,\cM_j,\lambda_j)$.  If the $\cM_j$
intersect correctly, and if for every bounded open $U$ $$\sum_{j=1}^\infty{\bf
M}[T_j\contract U]<\infty,\eqno(6.2)$$
then $\sum T_j$ is a weakly laminar current.  If the $T_j$ are laminar with
pairwise disjoint carriers, and if (6.2) holds, then $\sum T_j$ is laminar.
 
\give Proof.  We let $\cM$ (resp.\ $A$) denote the disjoint union of the
$\cM_j$ (resp.\ $A_j$), and we define the measure $\lambda=\sum\lambda_j$ by
setting $\lambda\contract{A_j}=\lambda_j$.  By (6.2), it follows that
(6.1) holds, so $(A,\cM,\lambda)$ represents a positive current, which must
coincide with $\sum T_j$.\qed

\give Example.  Weakly laminar currents are well behaved with respect to taking
summations, but for our applications we will need to take the supremum
of an increasing
family of laminar currents.  To understand some of the technical points of the
sequel, it may be helpful to note that although $T_1$ and $T_2$ are
uniformly laminar currents, and $T_1\le T_2$, it may happen that the positive
current $T_2-T_1$ is not weakly laminar.  Similarly,  $T_1+T_2$
and $\max(T_2,2T_1)$ may fail to be laminar.  For a simple example, consider
$T_1=[M_1]\le T_2=[M_2]$, where $M_1\subset M_2\subset\cx{}$, but $M_2\cap
\partial M_1$ has positive area. \medskip
 
Let us discuss the polar representation $T=t\nu$ of a laminar current.  From
(5.1) we have $[M]=t_M\cH^2\contract M$.  Thus the underlying measure is
$$\nu=\int_{a\in A}\lambda(a)\cH^2\contract M_a=|T|,\eqno(6.3)$$
and the set $\bigcup_{a\in A}M_a$ carries full measure for $\nu$.   By (6.3),
$\nu(E)=0$ holds for a Borel set $E$ if and only if $Area(M_a\cap E)=0$ for
$\lambda$ a.e.~$a$.  Since the manifolds $M_a$ intersect correctly, it follows
that for $\nu$ a.e. $x\in\bigcup_{a\in A}M_a$, the 2-vector is
$t(x)=t_{M_a}(x)$. Thus $t$ is a simple 2 vector at $\nu$ a.e.~point.  In
other words, there are vectors $t_1$ and $t_2$ such that $t=t_1\wedge t_2$.
The field of 2-vectors $t$ and $\nu$ depend only on $T$ and are independent of
the representation used to define them.

We let $\xi$ denote a family
of 1-dimensional complex manifolds $\alpha\subset\cx 2$ such that each
$\alpha\in\xi$ defines a current of integration $[\alpha]$ with finite mass
norm.  We will say that $\xi$ is a {\it stratified carrier} for a weakly
laminar current $T$ if \item{(i)}  $E:=\bigcup_{\alpha\in\xi}\alpha$ is a Borel
set. \item{(ii)}
$\xi$ is a measurable partition of $E$.
\item{(iii)} For $\lambda$ a.e.\ $M\in\cM$
there is a countable family $\{\alpha_i\}\subset\xi$ such that
$M-\bigcup_i\alpha_i$ has zero area.
 
If $T$ is laminar, then $\cM$ is a stratified carrier. 
It is a consequence of (6.3) that if $\xi$ satisfies (i), (ii), and (iii), then
$E$ is a carrier for $T$.  In Corollary 6.7 it will be shown that condition
(iii) is in fact independent of the choice of representation $(A,\cM,\lambda)$.
We note that the main difference between $\cM$ and $\xi$ is that the complex
manifolds in $\xi$ are disjoint.  We will say that two stratified carriers {\it
intersect correctly} if the complex manifolds in the stratifications intersect
correctly.  We say that a representation $(A,\cM,\lambda)$ is {\it subordinate
to} $\xi$ if for $\lambda$ a.e.\ $a\in A$ there exists $\alpha\in\xi$ with
$M_a\subset\alpha$.
 
The point of considering a stratified carrier is as follows.  Let us suppose
that a laminar current $T$ has a representation $(A,\cM,\lambda)$ which is
subordinate to a stratified carrier $\xi$.  (It
will be shown in Lemma 6.8 that any representation may be refined to be
subordinate to a given stratified carrier $\xi$.)  For $\alpha\in\xi$ we set
$A_\alpha=\{a\in A:M_a\subset \alpha\}$.  We
may let $\lambda_\xi$ denote the measure $\lambda$ restricted to the (coarser)
$\sigma$-algebra which is generated by $\xi$.  For $\lambda_\xi$\
almost every $\alpha\in\xi$ there is a conditional measure
$\lambda(\cdot|\alpha)$ on $A_\alpha$, as in \S2.  Let us define a function on
$\alpha$ by setting
$$\varphi^\alpha:=\int_{a\in A_\alpha}\chi_{M_a}\lambda(a|\alpha),\eqno(6.4)$$
where
$\chi_{M_a}$ denotes the function which is 1 on the set $M_a$ and 0 on
$\alpha-M_a$.  Since $M_a$ is an open subset of $\alpha$, and since the
conditional measure is positive, $\varphi^\alpha$ is lower semicontinuous on
$\alpha$. It is immediate that
$$\varphi^\alpha[\alpha]=\int_{a\in A_\alpha}[M_a]\lambda(a|\alpha).$$
It follows
from the defining property of the conditional measures that
$$T=\int_{\alpha\in\xi}\varphi^\alpha[\alpha]\lambda_\xi(\alpha).\eqno(6.5)$$
This differs
from the original representation of $T$ as a direct integral in that
the currents involved are not locally closed, but it has the advantage that the
supports may be taken to be essentially disjoint.

\proclaim Proposition 6.2. Let $T$ be a weakly laminar current.  Then for
$\epsilon>0$ and any bounded, open set $U$, there exist uniformly laminar
currents $T_j$ with disjoint supports such that
$${\bf M}\left[\left(T-\sum T_j\right)\contract U\right]<\epsilon.\eqno(6.6)$$
\vskip0pt
If $T$ is a laminar current, then there exist uniformly laminar currents
$T_1,T_2,\dots$ with disjoint supports such that $T=\sum T_j$.  Further, there
is a compact $K\subset U$ such that ${\bf M}[T\contract(U-K)]<\epsilon$ and
$T\contract K$ is the finite sum of uniformly laminar currents with disjoint
supports.
 
\present{Proof.} Let $T$ have a representation $(A,\cM,\lambda)$.  Let
$\cQ_n$ denote the decomposition of $\cx{}$ into squares of side $2^{-n}$ and
vertices at the points $(j+ik)2^{-n}$ for $j,k\in{\bf Z}$.  Let $\pi(x,y)=x$.
We may assume that the set of $a\in A$ such that $\pi(M_a)$ is a point has
$\lambda$ measure zero.  For each $a\in A$, we call a component $M^\prime$ of
$M_a\cap\pi^{-1}Q$ {\it good} if $\pi|_{M^\prime}:M^\prime\to Q$ is a
homeomorphism.  We let $\hat M_a(Q)$ be the union of all of the good components
 of
$M_a\cap\pi^{-1}Q$, and we set
$$T_Q=\int_{a\in A}\lambda(a)\,[\hat M_a(Q)].\eqno(6.7)$$
It is immediate that $$\sum_{Q\in\cQ_n}T_Q\le T.$$

Let $\cN$ denote the set of every  disk which arises as a good component of
$M_a\cap\pi^{-1}Q$ for some $a\in A$.  Thus there is a measure $\lambda_Q$ on
$\cN$ such that
$$T_Q=\int_{\cN}\lambda_Q(N)\,[N].$$
We observe that if $N_1,N_2\in\cN$, then the condition of correct intersection
implies that either $N_1\cap N_2=\emptyset$, or $N_1=N_2$.  Thus each $T_Q$ is
uniformly laminar.
 
We let
$$T^{(1)}=\sum_{Q\in\cQ_1}T^{}_Q,$$
so that $T^{(1)}$ is the sum of uniformly laminar currents with disjoint
carriers.
 
Now we suppose that $T^{(j)}_Q$ have been constructed for $1\le j\le n-1$ and
$Q\in \cQ_j$.  Each $T^{(j)}_Q$ is uniformly laminar, and $T^{(j)}=\sum_Q
T^{(j)}_Q$ is laminar.  Further, $T^{(1)}+\dots+T^{(n-1)}\le T$.  Since
$T-T^{(1)}-\dots-T^{(n-1)}$ is weakly laminar, we may let
$$T^{(n)}_Q:=(T-T^{(1)}-\dots-T^{(n-1)})_Q$$ be the uniformly laminar current
obtained in the construction (6.7).
 
We observe that if $U=\{|Re\,x|, |Re\,y|,
|Im\,x|, |Im\,y|<m\}$ for some integer $m$, then
 $(T^{(1)}+\dots+T^{(n)})\contract
U$ is a finite sum of uniformly laminar currents with disjoint carriers.  By the
construction above, the mass norm in (6.6) is given by
$$\int\lambda(a)\,Area\left(M_a\cap U-\bigcup_{Q\in\cQ_n}\hat
M_a(A)\right).$$
For fixed $a\in A$, the area decreases to zero as $n\to\infty$,
so this integral tends to zero by monotone convergence.

If $Q\in\cQ_n$, then $(T^{(1)}+\dots+T^{(n)})\contract(Q\times\cx{})$ is
uniformly laminar.  Thus the currents in the family
$\{T_j\}:=\{(T^{(1)}+\dots+T^{(n)})\contract(Q\times\cx{}):Q\in\cQ_n\}$ are
uniformly laminar and have disjoint carriers and satisfy (6.6) for $n$
sufficiently large.  If $T_j$ is restricted to a smaller compact inside its
carrier, the supports of $\{T_j\}$ will be parwise disjoint.

Finally, let us observe that if $T$ is laminar, then the carriers of
$T^{(j)}_Q$ are already pairwise disjoint, and
$T=\sum_{Q,j}T_Q^{(j)}$.  By subdividing the support of each $T_Q^{(j)}$ into
countably many compact sets, we have the first assertion.  The existence of $K$
with the required properties is a property of Radon measures. \qed

\give Remark.  It follows that the currents $\tilde T_1$ and $\tilde T_2$
defined in (5.7--8) are not locally weakly laminar on any open set.

\proclaim Lemma 6.3.  If $T_1,\dots,T_k$ are laminar currents with
representations that intersect correctly, then there exists $\xi$ which is a
stratified carrier for $T_j$ for $1\le j\le k$.
 
\give Proof.  The proof of this Lemma is a repetition of the proof of
Proposition 6.2 with $\cM$ replaced by $\cM_1\cup\dots\cup\cM_k$.  
To obtain a stratified carrier, we fix $n$ and $Q\in\cQ_n$.  We use the
notation $\cN^n_Q$ for the set $\cN$ defined above: the union over $a\in A$ of
the set of disks which are good components of $M_a\cap\pi^{-1}Q$.  We let
$\xi_1=\bigcup_{Q\in\cQ_1}\cN^1_Q$.  We continue inductively, setting
$\xi_n=\bigcup_{Q\in\cQ^n}\cN^n_Q-\xi_{n-1}$.  Finally, $\xi=\bigcup\xi_n$
has the desired properties.
\qed

Given a representation $(A,\cM,\lambda)$ of $T$, we may define a family of
germs of complex manifolds as follows: for $x\in\bigcup_{a\in A}M_a$, we let
$\hat M(x)$ be the germ of $x$ of the manifold $M_{a(x)}$ containing $x$.  The
correspondence $x\mapsto \hat M(x)$ is thus well defined $\nu$ a.e.~in terms of
the representation.  By (5.1) and (6.3), we have
$$T=\int_{a\in A}\lambda(a)\,t_{M_a}\sigma_{M_a}.$$
Since the $M_a$ overlap correctly, it follows
that if $(A^\prime,\cM^\prime,\lambda^\prime)$ and
$(A^{\prime\prime},\cM^{\prime\prime},\lambda^{\prime\prime})$ are two
representations, then
$$T_x\hat
M^\prime(x)=T_x\hat M^{\prime\prime}(x)\eqno(6.8)$$
holds for $\nu$ a.e.~$x$ (so
the germs intersect tangentially).  We now show that these germs coincide at
$\nu$ a.e.~point.
 
First we need a lemma.
 
\proclaim Lemma 6.4.  Let $M_1$ and $M_2$ be complex submanifolds of $\cx2$ such
that $M_1\cap M_2=\{p\}$, and $T_pM_1=T_pM_2$.  If $M_1^\prime$ is
sufficiently close to $M_1$, but $M_1^\prime\cap M_1=\emptyset$, then the
intersection $M_1^\prime\cap M_2$ is nonempty, and nontangential at all
intersection points.
 
\give Proof.  Let $k$ be the multiplicity of the intersection of $M_1$ and $M_2$
at $p$.  By the continuity of the intersection of complex manifolds, the
intersection of $M_1^\prime$ and $M_2$ (with multiplitity) near $p$ is $k$.
Thus it suffices to show that $M^\prime_1\cap M_2$ contains $k$ distinct points
near $p$.
 
Without loss of generality, we may work in a small neighborhood of $p=0$ and
assume that $\{y=f(x):|x|<1+\epsilon\}\subset\subset M_2$ for some holomorphic
function $f(x)=x^k+\dots\,$ and $\{|x|<1+\epsilon,y=0\}\subset\subset M_1$.  We
may assume that $\{y=f(x)=0:|x|<1+\epsilon\}=\{0\}$.  A manifold $M_1^\prime$
which is $C^1$ close to $M_1$ is of the form
$\{y=g(x):|x|<1+\epsilon\}\subset\subset M^\prime_1$.  The hypothesis that
$M_1\cap M_1^\prime=\emptyset$ implies that $g\ne0$.  By the Harnack
inequalities there is a constant $C_\epsilon$ such that
$$C_\epsilon^{-1}|g(0)|\le|g(x)|\le C_\epsilon|g(0)|$$
for $|x|\le 1$.  This implies that the higher order terms in $g(0)/g(x)=1+\dots$
are uniformly small.   Since $M_1^\prime\cap\{y=f(x):|x|<1\}$
is given by  $$g(0)=f(x){g(0)\over g(x)}=x^k+\dots,$$
and the higher order terms are uniformly bounded, this equation has $k$ distinct
solutions near $x=0$ for $g(0)$ sufficiently small.
\qed
 
\proclaim Lemma 6.5.  Let $(A^\prime,\cM^\prime,\lambda^\prime)$ and
$(A^{\prime\prime},\cM^{\prime\prime},\lambda^{\prime\prime})$ be two
representations for the weakly laminar current $T$.  Then $\hat M^\prime(x)=\hat
M^{\prime\prime}(x)$ for $\nu$ a.e. $x$.
 
\give Proof. Let $B=\{x:\hat M^\prime(x)\ne\hat
M^{\prime\prime}(x)\}$. Removing a set of
measure zero, we may assume that (6.8) holds at every point of $B$.   We must
 show
that $\nu(B)=0$. Otherwise, we may choose $\epsilon$ such that
$0<\epsilon<\nu(B)$ and let $T=\sum T_j^\prime$ be the sum of uniformly laminar
currents obtained in Lemma 6.2 corresponding to $\cM^\prime$.  If $T_j\contract
B=0$ for all $j$, then
$${\bf M}\left[T-\sum T_j^\prime\right]\ge{\bf M}[T\contract B]=\nu(B)$$
so it follows that  $T_j^\prime\contract B\ne0$ for some $j$.  Now the current
$T^\prime:=T^\prime_j$ is uniformly laminar and has the form
$$T^\prime=\int_{a^\prime\in
A^\prime}\lambda^\prime(a^\prime)[\Gamma^\prime_{a^\prime}].\eqno(6.9)$$ Let us
set $$B^\prime=\bigcup_{a^\prime\in A^\prime}\Gamma^\prime_{a^\prime}\cap B.$$
Since $T^\prime\contract B\ne0$, we have
$$|T^\prime|(B)=\int_{A^\prime}\lambda^\prime(a^\prime)
\,Area(\Gamma^\prime_{a^\prime}\cap B)>0,\eqno(6.10)$$ as in (6.3).  It follows
that $Area(\Gamma^\prime_{a^\prime}\cap B)>0$ for a set of positive
$\lambda^\prime$ measure, so $\nu(B^\prime)>0$. Now we let $T=\sum
T_j^{\prime\prime}$ be as in Lemma 6.2 for $\epsilon<\nu(B^\prime)$.  As before
there exists $k$ such that
$T^{\prime\prime}_k\contract B^\prime\ne0$.  Now we set
$T^{\prime\prime}:=T_k^{\prime\prime}$, and we represent $T^{\prime\prime}$ in a
form analogous to (6.9).  By the analogue of (6.10), we know that there exists
$a^{\prime\prime}$ such that $$Area(\Gamma^{\prime\prime}_{a^{\prime\prime}}\cap
B^\prime)>0.$$
 
Now let $b\in\Gamma^{\prime\prime}_{a^{\prime\prime}}\cap B^\prime$ be a point
of density with respect to area measure.  thus there is a sequence
$\{b_j\}\subset \Gamma^{\prime\prime}_{a^{\prime\prime}}\cap B^\prime$
converging to $b$.  Since $b,b_j\in B^\prime$, there exist
$a^\prime,a_j^\prime\in A^\prime$ with $b\in\Gamma^\prime_{a^\prime}$ and
$b_j\in\Gamma^\prime_{a_j^\prime}$.
Let $M_1=\Gamma^\prime_{a^\prime}$, $M_1^\prime=\Gamma^\prime_{a_j^\prime}$,
and $M_2=\Gamma^{\prime\prime}_{a^{\prime\prime}}$.  Since $b\in B$, $M_1$ and
$M_2$ define different germs of complex manifolds, and we may intersect them
with $B(b,\epsilon)$, if necessary, to have $M_1\cap M_2=\{b\}$.  Since (6.8)
holds at $b$, $M_1$, $M_1^\prime$, and $M_2$ satisfy the hypotheses of Lemma
6.3, so we conclude that all intersection points of $M_1^\prime$ and $M_2$ are
transversal.  But $b_j\in M^\prime_1\cap M_2$, and the intersection at $b_j$ is
tangential by (6.8).  By this contradiction we conclude that $\nu(B)=0$.
\qed

\proclaim Corollary 6.6.  Let $(\cM_1,\lambda_1)$ and $(\cM_2,\lambda_2)$ be
two representations of a weakly laminar current $T$.  Then for $\lambda_1$
a.e.~$M_1\in\cM_1$ and $\lambda_2$ a.e. $M_2\in\cM_2$, $M_1$ and $M_2$
intersect correctly.
 
In other words, the set of manifolds $\cM$ associated with a weakly laminar
current $T$ are unique, up to subdivision or refinement.
We say that weakly laminar
currents $T^\prime$ and $T^{\prime\prime}$ {\it intersect correctly} if they
 have
representations $(A^\prime,\cM^\prime,\lambda^\prime)$ and
$(A^{\prime\prime},\cM^{\prime\prime},\lambda^{\prime\prime})$ such that the
disks of $\cM^\prime$ and $\cM^{\prime\prime}$ intersect correctly.  By
 Corollary
6.6, this condition is independent of the representations $\cM^\prime$ and
$\cM^{\prime\prime}$ chosen.  Another consequence is the following.

\proclaim Corollary 6.7.  If $\xi$ is a stratified carrier which satisfies
condition (iii) for $\cM_1$, then (iii) holds for any other representation
$\cM_2$.
 
Let $(A,\cM,\lambda)$ be a representation of a weakly laminar current, and let
$\xi$ be a stratified carrier.  We will show how to subdivide the elements of
$\cM$ so that the representation is subordinate to $\xi$.  We set $\tilde A=
A\times\xi$, and we define a measure $\tilde \lambda$ on $\tilde A$ by setting
$\tilde\lambda(E\times\{\alpha\})=\lambda(E)$ for any measurable $E\subset A$
and any $\alpha\in\xi$.  In other words, $\tilde\lambda=\lambda\times\cH^0$ is
the product measure obtained from $\lambda$ and the counting measure $\cH^0$
on $\xi$.  We define $\tilde\cM$ by setting $\tilde
M_{(a,\alpha)}=M_a\cap\alpha$ for any $a\in A$ and $\alpha\in\xi$.

\proclaim Lemma 6.8.  If $(A,\cM,\lambda)$ is a representation of $T$, and if
$\xi$ is a stratified carrier of $T$, then $(\tilde A,\tilde\cM,\tilde\lambda)$
is a representation of $T$ which is subordinate to $\xi$.
 
\give Proof.  By definition of $\tilde A$, $\tilde\cM$, and $\tilde\lambda$, we
have
$$\eqalign{
\int_{\tilde a\in\tilde A}[\tilde M_{\tilde a}]\tilde\lambda(\tilde a) &=
\int_{(a,\alpha)\in A\times\xi}[M_a\cap\alpha]\tilde\lambda(a,\alpha)\cr
&=\int_{a\in A}\lambda(a)\int_{\alpha\in\xi}[M_a\cap\alpha]\cH^0(\alpha)\cr
&=\int_{a\in A}\lambda(a)\sum_{\alpha\in\xi}[M_a\cap\alpha]\cr
&=\int_{a\in A}\lambda(a)[M_a]=T\cr}$$
where the second line follows by the Fubini Theorem, and the fourth line is by
(iii) of the definition of stratified carrier. \qed
 
Since we may subdivide any representation $(A,\cM,\lambda)$ to be subordinate
to a given stratified carrier $\xi$, it follows that $T$ may be given as a
direct integral over the elements of $\xi$, as in (6.4) and (6.5).  This yields
the following:
 
\proclaim Lemma 6.9.  If $T$ is weakly laminar, and if $\xi$ is any stratified
carrier, then $T$ may be represented in terms of $\xi$ as follows: For
$\alpha\in\xi$, there exists a lower semicontinuous function
$\varphi^\alpha\ge0$ on $\alpha$ such that $\xi\ni\alpha\mapsto\varphi^\alpha$
is measurable, and
$$T=\int_{\alpha\in\xi}\lambda_\xi(\alpha)\,\varphi^\alpha[\alpha].$$
$T$ is laminar if and only if $\varphi^\alpha$ is locally constant a.e.\ on
$\{\varphi^\alpha>0\}$.

The {\it maximum}, written $\max(T_1,\dots,T_n)$, of the currents
 $T_1,\dots,T_n$
(if it exists) is characterized by the properties: $T_j\le\max(T_1,\dots,T_n)$
for $j=1,\dots,n$, and if $S$ is any current satisfying $T_j\le S$,
$j=1,\dots,n$, then $\max(T_1,\dots,T_n)\le S$.
 
\proclaim Lemma 6.10.  Let $T_1,\dots,T_n$ be weakly laminar currents which
intersect correctly.  Then $\max(T_1,\dots,T_n)$ exists as a positive current and
is weakly laminar.
 
\give Proof.  By Lemmas 6.3 and 6.8,  we may assume that the
representations of $T_j$ are subordinate to some carrier $\xi$. By Lemma 6.9
$$T_j=\int_{\alpha\in\xi}\varphi^\alpha_j[\alpha]\,\lambda^j_\xi(\alpha)$$ with
the measurable family of lower semicontinuous functions $\varphi^\alpha_j$ on
$\alpha$ being given by (6.4).  Let us define
$\lambda:=\lambda^1_\xi+\dots+\lambda^n_\xi$, and let $h_j$ be a measurable
function such that $\lambda^j_\xi=h_j\lambda$.  It follows
$$\int_{\alpha\in\xi}\max(h_1\varphi^\alpha_1,\dots,h_n\varphi^\alpha_n)[\alpha]
 \,\lambda(\alpha),$$
defines a laminar current which has the properties of $\max(T_1,\dots,T_n)$.
 \qed
 
\proclaim Lemma 6.11.  Let $T_1,\dots,T_n$ be uniformly laminar currents which
intersect correctly.  Suppose that for any $M_i\in\cM_i$ and $M_j\in\cM_j$,
$M_i\cap \partial M_j$ has zero area in $M_i$. Then $\max(T_1,\dots,T_n)$ exists
as a positive current and is laminar.
 
\give Proof. The existence of $\max(T_1,\dots,T_n)$ follows from Lemma 6.10. 
Since the relative boundaries have zero area, this current is laminar by Lemma
6.9. \qed
 
\proclaim Lemma 6.12.  Let $T_1\le T_2\le\dots$ be an increasing sequence of
weakly laminar currents whose mass is locally bounded.  Suppose that there
 exists
$\xi$ which is a stratified carrier for all $T_n$.  Then $\sup_n T_n$ exists as
a positive current and is weakly laminar.
 
\give Proof.  Each current $T_n$ may be written as
$$T_n=\int\varphi^\alpha_n[\alpha]\,\lambda^n_\xi(\alpha).$$
There exists a sequence of functions $g_n>0$ on $\xi$ such that $m=\sum
g_n\lambda^n_\xi$ is a probability measure.  Clearly $\lambda^n_\xi\ll m$ for
each $n$, so there exist measurable functions $h_n$ such that $\lambda^n_\xi=h_n
m$.  Further, since the currents $T_n$ are increasing, the functions
$\varphi^\alpha_n h_n$ are increasing in $n$ for fixed $\alpha$.  Thus the
function $$\tilde\varphi^\alpha:=\lim_{n\to\infty}\varphi^\alpha_n h_n$$ is
finite  for $m$ a.e.\ $\alpha$ (since the $T_n$ have locally bounded mass) and
 is
thus lower semicontinuous.  We conclude, then, that
$$T:=\int_{\alpha\in\xi}\tilde\varphi^\alpha[\alpha]\,m(\alpha)$$ is a geometric
current, which clearly has the property of $\sup_n T_n$. \qed
 
\give Remark.  Some of the properties of weakly laminar currents may be
summarized as follows.  Let $T$ be weakly laminar, and let $\cS(T)$ denote the
set of weakly laminar currents $0\le S\le T$.  Then the subset of $\cS(T)$
consisting of finite sums of uniformly laminar currents with disjoint supports
is dense in the local mass norm (Proposition 6.2).  If $0\le\psi\le 1$ is lower
semicontinuous, then $\psi\cS(T)\subset\cS(T)$ (Lemma 6.9).  Finally, $\cS(T)$
is convex and closed under countable maxima (Lemmas 6.3 and 6.12).

\bigbreak
\section 7. $\mu^+$ is a Laminar Current
 
 In this section, we show that $\mu^+$ is laminar.\footnote*{We wish to thank
Cliff Earle for telling us about the Ahlfors Covering Theorem, which is the
principal tool in the proof.}
 Let $\tilde D$ denote a 1-dimensional complex submanifold of $\cx2$, and let
 $D\subset\tilde D$
be a relatively compact domain with smooth boundary.  Let us suppose that
$$\mu^-|_{\tilde D}(D)=c>0{\rm\ and\ }\mu^-|_{\tilde D}(\bar D-D)=0.$$
It follows by [BS3], then, that
$$\lim_{n\to\infty}d^{-n}f^{*n}[D]=c\mu^+.\eqno(7.1)$$
Further, by general properties of the filtration (see [BS1, \S2]), we may choose
$R<\infty$ such that for all $n<\infty$
$$f^{-n}\tilde D\subset\{|y|<R\}\cup\{|y|<|x|\}.\eqno(7.2)$$
Let $Q\subset\cx{}$ be a connected open set.  For each $n$, we consider the
connected components $M$ of $(f^{-n}D)\cap(Q\times\cx{})$ and the
preimage components $D^\prime$ in the domain $D\cap f^n(Q\times\cx{})$.
If a component $D^\prime$ of $D$ is relatively compact in $D$, we say that
$D^\prime$ is an {\it island}; otherwise, it is a {\it tongue}.  Let
$\pi:\cx2\to\cx{}$ be the projection $\pi(x,y)=x$.  If $D^\prime$ is an
island, we say that it is a {\it good island} if the projection $\pi\circ
f^{-n}$ is univalent on $D^\prime$.
 
We let $\cG_n(Q)$ denote the set of components $M$ of $(f^{-n}D)\cap
Q\times\cx{}$ which are graphs over $Q$.  This corresponds to the set of
good islands, and each good island may be identified with the graph of
analytic function $\varphi:Q\to\cx{}$.  If we fix a point $x_Q\in Q$, then each
element of $\cG_n(Q)$ is uniquely determined by the value $\varphi(x_Q)$, i.e.
the intersection $M\cap(\{x_Q\}\times\cx{})$.
 
By (7.2), the union $\bigcup_n\cG_n(Q)$ is a normal family, and we let
$\cG(Q)$ consist of all graphs $\{y=\varphi(x):x\in Q\}$ which are obtained
as limits of sequences $\varphi_{n}\in\cG_n(Q)$.  Since $f$ is a
diffeomorphism, the components of $\cG_n(Q)$ are disjoint.  It follows from
the Hurwitz Theorem, then, that any two different graphs in $\cG(Q)$ are in
fact disjoint.  We let $A_Q\subset\cx{}$ denote the closed set of points
$\{\varphi(x_Q):\varphi\in\cG(Q)\}$.  For $a\in A_Q$ we let $M_j(a)$ denote
the element of $\cG(Q)$ passing through $(x_Q,a)$.
 
For each $n$ we define a measure $\lambda_{n}^+=d^{-n}\sum\delta_p$,
where the summation is taken over all $p\in\{x_Q\}\times\cx{}$ which are
parameters of elements $\varphi\in\cG_n(Q)$.  For each $Q$, we choose a
subsequence $\{n_k\}$ such that the limit
$\lim_{k\to\infty}\lambda^+_{n_k}$
exists.  We let $\lambda^+_{Q}$ denote this limit, and it follows that
 $\lambda^+_{Q}$
supported on $A_Q$.  Now we define
$$\mu^+_{Q}=c^{-1}\int_{a\in A_Q}\lambda^+_{Q}(a)\,[M(a)].\eqno(7.3)$$
 
It is evident that
$$d^{-n}[f^{-n}D]\contract(Q\times\cx{})\ge d^{-n}\sum_{M\in\cG_n(Q)}
[M].\eqno(7.4)$$ Thus, passing to the limit through the subsequence
$\{n_k\}$, we have $$c\mu^+\contract(Q\times\cx{})\ge c\mu^+_{Q}.$$

Let us use the following notation.   For $k\ge0$ we let $\cQ_k$ denote a
dyadic subdivision of the complex plane $\cx{}$ into open squares with vertices
of the form $r2^{-k}+is2^{-k}$ with $r$ and $s$ both odd.  Let $\cQ_k^{(2)}$,
$\cQ_k^{(3)}$,  and $\cQ_k^{(4)}$ denote the three different translates of
$\cQ_k^{(1)}$, so that $\cQ_k=\bigcup_{\sigma=1}^4\cQ_k^{(\sigma)}$.
 
As before, we construct families of graphs $\cG_k(Q)$ for
$Q\in\cQ_k^{(\sigma)}$, for each $\sigma=1,2,3,4$.  If we write
$$\mu^+_k=\sum_{Q\in\cQ_k}\mu^+_Q,$$
then it is evident from (7.4) and (7.3) that
$$\mu^+_1\le\mu^+_2\le\dots\le\mu^+_k\le\dots\le\mu^+.\eqno(7.5)$$
 
Now suppose that $j>k$, $Q\in\cQ$, $Q^\prime\in\cQ_k$, and $Q\subset
Q^\prime$.  If $\cG(Q^\prime)|_{Q}$ denotes the restriction of the disks to
$Q$, then it is evident that $\cG(Q^\prime)|_{Q}\subset\cG(Q)$. Similarly,
making the natural identification via the holonomy for the transversal
measures, it follows that $\lambda^+_{Q^\prime}\le\lambda^+_{Q}$.  Thus if we
set
$$\tilde\lambda^+_{Q}:=\lambda^+_Q-\max\{\lambda^+_{Q^\prime}:Q^\prime\supset
Q,Q^\prime\ne Q\},$$ then $\tilde\lambda^+_Q$ is a positive measure.  For each
$Q$, then, we set $$\tilde\eta^+_Q=\int_{a\in A_Q}\tilde\lambda^+_Q(a)
[M_Q(a)]$$ and $$\tilde\eta^+_j=\sum_{Q\in\cQ}\tilde\eta^+_Q.$$
Thus by Lemma 6.1 we have shown:
 
\proclaim Lemma 7.1.  The currents $\tilde\eta^+_j$ are uniformly laminar over
the squares of $\cQ_j$ and have disjoint carriers, and
$$\lim_{k\to\infty}\mu^+_k=\sum_{j=1}^\infty\eta^+_j.$$
Further, this limit is a laminar current.
 
In Theorem 7.4 we will  show that this limit is equal to $\mu^+$.
 
Let $Q_1,\dots,Q_q\subset\cx{}$ be simply connected, open sets such that
$\overline Q_i\cap\overline A_j=\emptyset$ for $i\ne j$.  We let
$Q:=Q_1\cup\dots\cup Q_q$, and $I_0:=Area(Q)$.  We set $D_{(n)}:=D\cap
f^n(Q\times\cx{})$, and we consider the map $$g_n:=\pi\circ
f^{-n}:D_{(n)}\to\cx{}.$$ We let $I_{(n)}$ denote the area (with multiplicity)
of $g_n(D_{(n)})$.  The {\it mean sheeting number} of the map $g_n$ is
$S_{(n)}:=I_{(n)}/I_0$.  The length of the {\it relative boundary} is defined by
$$L_{(n)}:=Length(g_n(\partial D)\cap Q).$$ Fixing the number $n$ of iterates,
we write $N(Q_j)$ for the number of good islands over $Q_j$, i.e.~this is just
the cardinality of the set $\cG_n(Q_j)$.  We will use the following celebrated
result of Ahlfors (see Nevanlinna [N, Chapt.~XIII], or Hayman [Ha]).
 
\proclaim Ahlfors' Covering Theorem.  There is a constant $h$ depending only on
$Q$ such that the mappings $g_n$, $n=1,2,\dots$ satisfy
$$\sum_{j=1}^qN(Q_j)\ge(q-2)S_{(n)}-hL_{(n)}.$$
 
We will use this inequality to estimate the amount of mass in
$\sum\mu^+_{Q_j}$.  By [BS1] we have
$$\lim_{n\to\infty}d^{-n}I_{(n)}=c\,Area(Q),{\rm\ or\
}\lim_{n\to\infty}d^{-n}S_{(n)}=c,\eqno(7.6)$$
with $c$ as in  (7.1).  Further, by [BS3], there is a constant $C<\infty$ such
that
$$L^2_{(n)}\le Cd^n.\eqno(7.7)$$
 
We note that for a current $T$, the mass norm of $T\contract{i\over2}dx\wedge
d\bar x$ on the set $B\times\cx{}$ is the same as ${\bf M}[\pi_{*}T\contract
B]$.   Each $M_j(a)\in\cG_n(Q_j)$ is the graph of an analytic function
on $Q_j$.  Thus the mass norm is
$${\bf
M}\left[\,[M_j(a)]\contract{i\over2}dx\wedge d\bar x\right]=Area(Q_j).$$

\proclaim Lemma 7.2.  If $Area(Q_1)=\dots=Area(Q_q)$, then
$${\bf M}\left[\sum_{j=1}^q\mu^+_{Q_j}\contract{i\over2}dx\wedge d\bar
x\right]\ge {q-2\over q}Area(Q).$$
 
\give Proof.
By the definition (7.3), it follows that the mass norm of $\mu^+_{Q_j}$ is
$${\bf M}[\mu^+_{Q_j}\contract{i\over2}dx\wedge d\bar x]=c^{-1}{\bf
M}[\lambda^+_j]Area(Q_j).\eqno(7.8)$$
 
In order to estimate ${\bf M}[\lambda^+_j]$, we count the number of components
$M_{}$ that appear in the right hand side of (7.4).  This is the same as
the number of good islands over $Q_j$.  Thus we have
$$\eqalign{
{\bf
M}\left[\sum_jd^{-n}\sum_{M_{}\in\cG_n(Q_j)}[M_{}]\contract{i\over2}dx
\wedge
d\bar x\right]&\ge d^{-n}\sum_{j=1}^q\#\cG_n(Q_j)Area(Q_j)\cr 
&\ge d^{-n}{Area(Q)\over q}\sum_{j=1}^qN(Q_j)\cr
&\ge {d^{-n}\over q}{Area(Q)}((q-2)S_{(n)}-hL_{(n)}),\cr }$$
where the middle inequality follows
from the identity $Area(Q_j)=q^{-1}Area(Q)$, and the last inequality follows
from the Ahlfors Covering Theorem.  Applying (7.7), we have
$${\bf
M}\left[\sum_jd^{-n}\sum_{M_{}\in\cG_n(Q_j)}[M_{}]\contract{i\over2}dx\wedge
d\bar x\right]\ge{q-2\over q}Area(Q)(d^{-n}S_{(n)}-O(d^{-{n\over 2}})).$$
Letting $n\to\infty$, we see from (7.6) that the right hand side tends to
$c(q-2)Area(Q)/q$ as the left hand side tends to $\sum_j{\bf
M}[\lambda^+_j]Area(Q_j)$.  Combinded with (7.8), this yields Lemma 7.2.\qed

\proclaim Lemma 7.3.  Let $B\subset\cx{}$ denote the unit square.  Then
$${\bf M}[\pi_{*}(\mu^+-\mu^+_k)\contract B]\le 8\cdot 4^{-k}.$$
 
\give Proof.  We note that ${\bf M}[\pi_{*}\mu^+\contract B]=Area(B)$ for any
open set.  And since $\mu^+\ge\mu^+_k$,
$${\bf M}[\pi_{*}(\mu^+-\mu^+_k)\contract B]={\bf M}[\pi_{*}\mu^+\contract
B]-{\bf M}[\pi_{*}\mu^+_k\contract B].$$
Thus the Lemma follows by setting $q=4^{k-1}$ and adding the estimate of Lemma
7.2 over the four partitions $\cQ^{(\sigma)}_k$.\qed
 
\proclaim Theorem 7.4.  $\lim_{k\to\infty}\mu^+_k=\mu^+$, and $\mu^+$ is a
laminar current.
 
\give Proof.  If we show that the limit holds, then $\mu^+$ is laminar by
Lemma 7.1.  By (7.5), it suffices to show that
$$\lim_{k\to\infty}{\bf M}[\mu^+_k\contract\pi^{-1}B_0]={\bf
M}[\mu^+\contract\pi^{-1}B_0]$$
for any open $B_0\subset\cx{}$.  Without loss of generality, we may choose
$B_0$ to be relatively compact in the unit square $B$.
 
For $\alpha\in\cx{}$,
we define the projection $\pi^\prime(x,y)=x-\alpha y$.  Let us choose
$\alpha\ne0$ sufficiently small that
$$(\pi^{-1}B_0)\cap{\rm spt}\mu^+\subset{\pi^\prime}^{-1}B.$$
Following the procedure for constructing the current $\mu^+_k$, except that the
projection $\pi^\prime$ is used in place of $\pi$,  we may construct a current
${\mu^\prime}^+_k$.  Thus we use the function $g_n:=\pi^\prime\circ f^{-n}$,
and $\cG^\prime_n(Q_j)$ consists of manifolds which are graphs with respect to
the coordinates $x^\prime=x-\alpha y$ and $y^\prime=y$.  Corresponding to Lemma
6.2, we have
$${\bf M}[\pi^\prime_*(\mu^+-{\mu^\prime}^+_k)\contract B]\le
8\cdot4^{-k}.$$
 
By Lemma 6.6 there is a geometric current $T_k$ such that $\mu^+_k$,
${\mu^\prime}^+_k\le T_k\le\mu^+$.  Thus we have
$${\bf M}[(\mu^+-T_k)\contract\chi_{\pi^{-1}(B_0)}{i\over2}d(x-\alpha y)\wedge
d\overline{(x-\alpha y)}]\le 8\cdot 4^{-k}.$$
Now we use  the values $\alpha=0$ and $\alpha=\pm a\in{\bf R}$ and the identity
$$d(x-ay)\wedge d\overline{(x-ay)} + d(x+ay)\wedge d\overline{(x+ay)}-2dx\wedge
d\bar x=2dy\wedge d\bar y$$
to obtain
$$|a|^2{\bf M}[(\mu^+-T_k)\contract\chi_{\pi^{-1}(B_0)}{i\over 2}dy\wedge
d\bar y]\le 16\cdot 4^{-k}.$$
Thus
$${\bf M}[(\mu^+-T_k)\contract\chi_{\pi^{-1}(B_0)}\beta]\le
8(1+2|a|^{-2})4^{-k}$$
where $\beta={i\over2}(dx\wedge d\bar x+dy\wedge d\bar y)$.  Since $\mu^+-T_k$
is positive, this gives
$${\bf M}[(\mu^+-T_k)\contract\chi_{\pi^{-1}(B_0)}]\le 8(1+2|a|^{-2})4^{-k}.$$
Thus $\lim_{k\to\infty}T_k=\mu^+$.
 
Now let us recall that $T_k$ is obtained by taking $\mu^+_k$ and adding all of
the currents of integration that appear in ${\mu^\prime}^+_k$, after removing
the sets where a manifold $M\in\cG(Q_k)$ overlaps a manifold $M^\prime\in
\cG^\prime(Q^\prime_k)$.  But let us consider such a manifold $M^\prime\in
\cG^\prime(Q^\prime_k)$.  As we increase $k$ to a larger index, say $K$, we
subdivide it into the pieces $\pi^{-1}(Q)\cap M^\prime$ for
$Q\in\cQ^{(\sigma)}_K$.  For any point $P\in M^\prime$, except at the (finite)
set where $\pi$ is branched, there is a square $Q\in\cQ^{(\sigma)}_K$ for some
large $K$ such that a component of $\pi^{-1}\bar Q\cap M^\prime$ contains $P$,
and this component belongs to $\cG(Q)$.  Thus it follows from monotone
convergence that $\lim_{j\to\infty}\mu^+_j\ge T_k$.  Thus
$\lim_{k\to\infty}\mu^+_k=\mu^+$.
\qed
 
Let us denote the total space of
the graphs in $\cG(G)$ as $\cE(G)=\bigcup_{\Gamma\in\cG(G)}\Gamma$.
We may write $\mu^+$ in the polar form $\mu^+=t|\mu^+|$, where
$|\mu^+|$ is the total variation measure, and for $|\mu^+|$ a.e.\ point $p$,
$t(p)$ is the unit 2-vector tangent to  $M(p)\in\cM$.  Thus we may define
$$m^+:=\mu^+\contract({i\over2}dx\wedge d\bar x)=\langle {i\over2}dx\wedge d\bar
x,t(p)\rangle |\mu^+|.$$
We note, further, that the integral of $\langle {i\over2}dx\wedge d\bar
x,t(p)\rangle$ over a complex manifold $M$ is just the area (with
multiplicity) of the projection of $M$ to the $x$-axis.  Since $\mu^+$ is
laminar, and since $\langle {i\over2}dx\wedge d\bar
x,t(p)\rangle$ does not vanish identically on any stable manifold, it follows
that $m^+$ and $|\mu^+|$ define the same measure class.

\proclaim Theorem 7.5.  Let $G_1$, $G_2$, and $G_3$ be Jordan
domains in $\cx{}$ with disjoint closures.  Then for some $j$,
$$m^+({\cE(G_j)})\ge{1\over 9}Area(G_j).$$

\give Proof.  Let us recall the current $\mu^+_{Q_j}$, constructed above.
The total variation measure associated with this current
satisfies $|\mu^+_{G_j}|\le|\mu^+|$ (with $G_j=Q_j$).  It follows that
$$m^+({\cE(G_j)}) = \int_{\cE(G_j)}\langle
t_\Gamma(p),{i\over 2}dx\wedge d\bar x\rangle |\mu^+|(p)\ge{\bf
M}[\mu^+_{Q_j}\contract{i\over 2}dx\wedge d\bar x].$$

Without loss of generality, we may enlarge $G_j$ to a larger
Jordan domain $G_j^\prime$, so that the three domains have the same area.  If
we set $q=3$ in Lemma 7.2, then we have
$$\sum_{j=1}^3 m^+({\cE(G_j)})\ge{1\over 3}Area(G^\prime),$$
where $Area(G^\prime)$ is the area of any of the $G_j^\prime$.
It follows, then, that for some $j$,
$$m^+({\cE(G_j^\prime)})\ge{1\over
9}Area(G_j^\prime).$$ Finally, since each $\Gamma^\prime$ is a graph over the
(larger) domain $G^\prime_j$, this inequality remains after we shrink to the
domain $G_j$.\qed
 
\proclaim  Theorem 7.6 (Three Islands). Let $G_1$, $G_2$, and $G_3$ be Jordan
domains with disjoint closures.  Then for some $j$, the total space of
$\cG(G_j)$ has positive $|\mu^+|$ measure.
 
\give Remark.  In \S8, Corollary 8.8, it will be shown that (almost every)
manifold making up the laminar structure of $\mu^+$ is in fact an open subset of
one of the stable manifolds $W^s(p)$, $p\in\cR$ given by the Pesin theory. The
utility of this theorem is that it gives the existence of stable manifolds that
are graphs over arbitrarily large sets.

A more general formulation is as follows.  Let $h$ be any polynomial, and let
$\cG(G,h)$ denote the set of all components $M$ of
manifolds obtained in the construction of $\mu^+$ such that $h|_M:M\to h(M)$ is a
conformal equivalence.  Thus, with our previous notation, we have
$\cG(G)=\cG(G,\pi_x)$. Thus we have:
 
\proclaim Corollary 7.7. Let $G_1$, $G_2$, and $G_3$ be Jordan
domains with disjoint closures, and let $h$ be any polynomial.  Then for some
$j$, the total space of $\cG(G_j,h)$ has positive $|\mu^+|$ measure.

\bigbreak
 
\section 8. Geometric intersection of $\mu^+$ and $\mu^-$
 
 By Sections 6 and 7 we know that $|\mu^\pm|$ almost every point lies inside a
uniformly laminar current which makes up part of $\mu^\pm$.  In this chapter we
will obtain a uniformly laminar structure for the currents $\mu^+$ and $\mu^-$
near any regular point for $\mu$. This is possible due to a  hyperbolic structure
given by Pesin boxes. Given a Pesin box $P$,
we can identify it with $P^s\times P^u$ via an appropriate homeomorphism
(see \S4).  Then
by Theorem 4.7, $\mu$ also has a product structure on this box,
i.e.\ $\mu\contract P$ is taken via this homeomorphism to
$\lambda^s\otimes\lambda^u$, where the measures $\lambda^s$ and $\lambda^u$
are induced by the currents $\mu^-$ and $\mu^+$ correspondingly.
Let us fix an ``origin" $o\in P$. For any $a\in P^u, b\in P^s$,
 denote by $\Gamma^s(a)$ a piece of
$W^s_{loc}(a,b)$ which is projected onto the disk $B^s(o,r)$
parallel to $E^u(o,r)$ (it does not depend on $b$).
 Similarly we can define a family of disks $\Gamma^u(b)$ .
Now let us consider the following sets supplied with a uniformly
laminar structure:
 $$\Gamma^s=\bigcup_{a\in P^u}\Gamma^s(a),\qquad
    \Gamma^u=\bigcup_{a\in P^s}\Gamma^u(a).$$
If a Pesin box $P_j$ is labeled by $j$, we will use the same label
for the corresponding sets $\Gamma^s_j(a)$ etc.
We let $\{  P_j,j=1,2,\dots\}$ be
a family of Pesin boxes such that $\bigcup P_j$ has full measure,   and we
set   $$\eta^+_j:=\int_{b\in P^u_j}\lambda^u_j(b)\,[\Gamma_j^s(b)]$$
$$\eta^-_j:=\int_{a\in P^s_j}\lambda^s_j(a)\,[\Gamma_j^u(a)],$$
 which are
uniformly laminar currents.  Without loss of generality, we may assume that
these currents satisfy the hypotheses of Lemma 6.11.  Thus the currents
$$\eta^\pm_{[n]}=\max(\eta_1^\pm,\dots,\eta^\pm_n)\quad {\rm\ and\ }\quad
\eta^\pm=\lim_{n\to\infty}\eta^\pm_{[n]}$$ exist and are laminar.  By the
holonomy invariance obtained in \S4, it follows that $\eta_j^\pm$ is well
defined independently of the transversal used in the definition.

\proclaim Lemma 8.1.  The sets $\Gamma^{s/u}_j$  satisfy
 $\eta_j^+=\mu^+\contract \Gamma^s_j$ and $\eta_j^-=\mu^-\contract \Gamma^u_j$.
 Thus $\eta^\pm_j\le\mu^\pm$.
 
\give Proof. Let $M$ be any transversal to the lamination $\Gamma^s_j$.
Since the measure $\lambda^u_j$ is induced by the current $\mu^+$,
$$\eta^+|_M\leq\mu^+|_ M. \eqno (8.1)$$
Hence  $\eta^+\le\mu^+$.
 
Let $\eta^+=\tau|\eta^+|$ and  $\mu^+=t|\mu^+|$ be the polar representations.
Then $|\eta^+|\le|\mu^+|$.  Since $\mu^+$ is a laminar current, $t$ is a
simple 2 vector $|\mu^+|$ a.e.  Thus $\tau=t$ $\; |\eta^+|$-a.e.,
and the Lemma follows. \qed

\proclaim Lemma 8.2.  If $T$ is a closed current, $0\le T\le \mu^+$, then
locally there is a continuous function $u$ with $dd^cu=T$.
 
\give Proof.  Since $T$ is closed, there is locally an integrable function $u$
such that $dd^cu=T$.  If $\beta={i\over2}(dx\wedge d\bar x+dy\wedge d\bar y)$,
then $\Delta u\contract\beta=dd^cu$.  It follows that $0\le\Delta u\le \Delta
G^+$.  Let $\nu$ denote the positive measure $\Delta G^+-\Delta u$ on some open
set $\cO$, and let $s=-c_4 |x|^{-2}*\nu$ denote the convolution with $\nu$, with
$c_4$ chosen so that of $-c_4|x|^{-2}$ is the fundamental solution of $\Delta$
on  ${\bf R}^4$.  Thus $s$ is
subharmonic, and the difference between $G^+$ and $u+s$ is harmonic on
$\cO\subset{\bf R}^4$.  A subharmonic function $v$ on $\cO$ satisfies
$\liminf_{q\to q_0}v(q)\le\limsup_{q\to q_0}v(q)=v(q_0)$ for all
$q_0\in\cO$.  Since $s$ and $u$ both satisfy this inequality, and since
$u+s$ is continuous at $q_0$, it follows that $s$ and $u$ are continuous at
$q_0$.\qed
 
We will define two different ways of taking the product of two currents.
First, we consider a continuous, psh function $u$ and a positive,
closed (1,1) current $T$.  We define the (2,2) current $T\wedge dd^cu$ by
its action on a test function $\varphi$:
$$(T\wedge dd^cu)(\varphi)=T(u\,dd^c\varphi).$$
(This is essentially just integrating the $dd^c$ by parts since $T$ is closed.)
It is evident from the right hand side of the defining equation that if $u_j$
converges uniformly to $u$, then $T\wedge dd^cu_j$ converges to $T\wedge
dd^cu$.  We refer the reader to [BT1] for further discussion of the $\wedge$
operation.
 
If $L_1$ and $L_2$ are uniformly laminar currents on $\Delta^2$, then it is also
 natural
to define
$$L_1\,\dot\wedge\, L_2=\int\lambda_1(a_1)\int\lambda_2(a_2)\,
[\Gamma_{a_1}\cap\Gamma_{a_2}]$$
with $[\Gamma_{a_1}\cap\Gamma_{a_2}]$ defined as the 0-current which puts unit
mass on each point of $\Gamma_{a_1}\cap\Gamma_{a_2}$, with the exception that
$[\Gamma_{a_1}\cap\Gamma_{a_2}]=0$ if $\Gamma_{a_1}=\Gamma_{a_2}$.  This is
analogous to the integrated version of (5.3), except that
$\Gamma_{a_1}\cap\Gamma_{a_2}$ is not necessarily transversal or finite.
 
\proclaim Lemma 8.3.  Let $L$ and $L^\prime$ be uniformly laminar currents on
$\Delta^2$ such that there is a continuous, psh function $u$ with
$dd^cu=L$.  Then $$L\wedge L^\prime=L\,\dot\wedge\, L^\prime.$$
 
\present{Proof.}  Without loss of generality, we may assume that $L$ and
$L^\prime$ are represented in the form (5.2), and
$$u={1\over2\pi}\int\lambda(a)\log|y-\varphi_a(x)|.$$
It will suffice to work over the relatively compact set $\{|x|<1-\epsilon\}$.
Let us fix $\Gamma^\prime=\Gamma_{a^\prime}$.  Choosing a parameter $\zeta=x$
for points $(x,y)=(\zeta,\varphi_a(\zeta))\in \Gamma^\prime$, we have
$$\log|y-\varphi_a(x)|=\sum_{j=1}^{N_a}\log|\zeta-p_j(z)| +h_a(\zeta),$$
where $h_a$ is harmonic.  Since $h_a$ is harmonic on $\{|x|<1\}$, it is bounded
on $\{|x|<1-\epsilon\}$.  Let us define
$$A_R=\{a\in A:\|h_a\|_{L^\infty(|x|<1-\epsilon)}\le R,N_a\le R\}.$$
If we set
$$u_R(x,y)={1\over2\pi}\int_{a\in A_R}\lambda(a)\log|y-\varphi_a(x)|,$$
then, as in Lemma 8.2, $u_R$ is continuous.  Further, since the $A_R$ increase
to $A$ as $R\to\infty$, $u_R$ converges uniformly to $u$.  Thus
$$\eqalign{
(dd^cu_R)\wedge[\Gamma^\prime]&=(dd^c\int\lambda(a)\log|y-\varphi_a(x)|
\wedge[\Gamma^\prime]\cr
&=\left(\int\lambda(a)\,[\Gamma_a]\right)\wedge[\Gamma^\prime]\cr
&=\int\lambda(a)\,[\Gamma_a\cap\Gamma^\prime]=dd^cu_R\,\dot\wedge\,[\Gamma^\prime],\cr}$$
where the next to last equality follows from the Fubini theorem, since the
multiplicity of the intersection is uniformly bounded for $a\in A_R$.
Letting $R\to\infty$, we have
$$L\wedge[\Gamma_{a^\prime}]=L\,\dot\wedge\,[\Gamma_{a^\prime}].\eqno(8.2)$$
 
Finally, we integrate (8.2) with respect to $\lambda^\prime(a^\prime)$.  The
right hand side yields $L\,\dot\wedge\, L^\prime$ by Fubini's Theorem.  The left
hand side, applied to a smooth test function $\chi$ is
$$\eqalign{
\int\lambda^\prime(a^\prime)(L\wedge[\Gamma_{a^\prime}])\chi
&=\int\lambda^\prime(a^\prime)[\Gamma_{a^\prime}]u\,dd^c\chi=\cr
&=L^\prime(u\,dd^c\chi)=(L^\prime\wedge dd^cu)(\chi),\cr}$$
which completes the proof.\qed

\proclaim Lemma 8.4. We have $\eta_j^+\wedge\eta_j^-=\mu\contract P_j$, and
thus $\mu^-\wedge\eta^+_{[k]}\ge\mu\contract\bigcup_{j=1}^k P_j.$
 
\give Proof.  By Lemmas 8.2 and 8.3, we have
$\eta^+_j\wedge\eta^-_j=\eta^+_j\,\dot\wedge\,\eta^-_j$.
 By the product structure of
Theorem 4.7, we have that under the homeomorphism between $ P_j$ and
$P_j^s\times P_j^u$, $\eta_j^+\,\dot\wedge\,\eta_j^-$ is taken to
$\lambda^s_j\otimes\lambda^u_j$, which in turn is equivalent to
$\mu\contract P_j$.  Similarly, since $\mu^-\ge\eta_j^-$ and
$\eta^+_{[k]}\ge\eta^+_j$, we have
$$\mu^-\wedge\eta^+_{[k]}\ge\eta_j^+\wedge\eta_j^-
=\eta^+_j\dot\wedge\eta^-_j=\mu\contract P_j.$$ Since this holds for all $j$,
the Lemma follows.\qed
 
\proclaim Lemma 8.5.  $\lim_{n\to\infty}d^{-n}f^{*n}\eta^+=\mu^+$.
 
\give Proof.  Let $\varphi$ be a test form.  We will show that
$$\int\varphi\mu^+=\lim_{n\to\infty}\int\varphi d^{-n}f^{*n}\eta^+.
\eqno(8.3)$$
Without loss of generality, we may assume that $\varphi\ge0$.  Let $\epsilon>0$
be given.  By Lemma 8.4, we may choose $k$ large enough that the total mass of
$\mu^-\wedge\eta^+_{[k]}$ is greater than $1-\epsilon$.  By Lemma 6.1, we may
write $\eta^+_{[k]}$ as a sum of uniformly laminar currents $\sum L_j$
with disjoint
carriers.  We may take finitely many terms from this summation and choose test
functions $0\le\chi_j\le1$ such that $\sum\chi_j L_j\le\eta^+_{[k]}$, and the
total mass of $\mu^-\wedge\sum\chi_jL_j$ is $c>1-\epsilon$.
 
Now by [BS3], we have
$\lim_{n\to\infty}d^{-n}f^{*n}(\sum\chi_jL_j)=c\mu^+$.  Since $\varphi\ge0$, we
have
$$\int\varphi\mu^+\ge\int\varphi d^{-n}f^{*n}\eta^+\ge$$
$$\ge\int\varphi
d^{-n}f^{*n}\left(\sum\chi_jL_j\right)\ge(1-\epsilon)\int\varphi\mu^+$$ for $n$
sufficiently large.  Since $\epsilon$ may be made arbitrarily small, we have
(8.2).\qed

Let us assume further that for $\lambda^s$
 a.e.\ $a$, the measure induced by $\mu^-$ on the corresponding stable
manifold puts
no mass on $\partial\Gamma_j^s(a)$. Then we have the following:
 
\proclaim Lemma 8.6.  Let $\tilde M$ be a 1-dimensional submanifold of
$\cx2$, and let $M$ be a relatively compact submanifold such that
$\mu^+|_{\tilde M}(\partial M)=0$.  Then
$$\lim_{n\to\infty}(d^{-n}f^{*n}\eta_j^+)\,\dot\wedge\,[M]=c\mu^+\wedge[M],$$
where $c:=\eta^+_j|_{\tilde M}[M]$.

 If we set $G^s=\cup_j\Gamma^s_j$, then by Lemma 8.1 we have
$$d^{-n}f^{*n}\eta^+=\mu^+\contract f^{-n}(G^s).$$
Since $G^s\subset\bigcup_{x\in\cR}W^s(x)$, where $\cR$ is the set
of all regular points (see \S2), it follows from  Lemma
8.5 that we have:
 
\proclaim Corollary 8.7. $\bigcup_{x\in\cR}W^s(x)$ is a
carrier for $|\mu^+|$.
 
By Corollary 6.6, we have:
 
\proclaim Corollary 8.8.  If $(A,\cM,\lambda)$ is a representation of $\mu^+$,
then $\lambda$ almost every $M\in\cM$ is an open subset of a stable manifold
$W^s(x),\; x\in \cR$.
 
Here we give a slightly different formulation of holonomy invariance.  This
is more general than the one given in \S4 because it applies to all stable
manifolds.  Let $\cM=\{M_\alpha:\alpha\in A\}$ be a family of stable manifolds.
Let $\cD=\{D_t:0\le t\le 1\}$ be a continuous family of manifolds such that each
$D_t$ is a transversal to $\cM$.  We define $X_j= D_j\cap\bigcup_\alpha
M_\alpha$, $j=1,2$.  The holonomy map $\chi:X_0\to X_1$ is defined by at a point
$x\in X_0$ by following the intersection point with $D_t$ from $t=0$ to $t=1$.
 
\proclaim Theorem 8.9 (Holonomy Invariance).  The holonomy map preserves the
slices of $\mu^+$, i.e. $$\chi_*(\mu^+|_{D_0}\contract
X_0)=\mu^+|_{D_1}\contract X_1.$$
 
\give Proof.   If the family
$\cM$ consists of leaves of $ G^s=\cup_j\Gamma^s_j $,
then holonomy is preserved.  In
general, we consider the compact sets $\gamma_\alpha=\{D_t\cap M_\alpha:0\le
 t\le
1\}$.  For each $\alpha$ the curve $\gamma_\alpha$ is contained in a stable
manifold, so there is an $n$ such that $f^n\gamma_\alpha$ is contained in one of
the leaves of $ G^s$.  Thus for $\epsilon>0$ there exists an $n$ such that
$\{x\in X_0:f^n\gamma_\alpha\not\subset G^s\}$ has measure less than
$\epsilon$.  Since the holonomy is preserved on the complement of this set, we
see that the Lemma holds.\qed

\bigbreak

\section 9. Saddle Points and the Support of $\mu$

 Questions about saddle points have motivated much of the preceding work on
currents, and we are grateful to J.H.\ Hubbard for several discussions on
this subject.  In this section we show that there is  a homoclinic/heteroclinic
 intersection
between any pair of stable and unstable manifolds.  The general idea is that if
$D^s$ is a stable disk through a saddle point $p$, then the normalized pullbacks
$d^{-n}[f^{-n}D^s]$ converge to a nonzero multiple of $\mu^+$.  By the results
of \S8, it follows that the product $d^{-n}[f^{-n}D^s]\wedge\mu^-$ 
converges to a nonzero multiple of the measure $\mu=\mu^+\wedge\mu^-$.  Since
this is also equal to the intersection wedge product, it follows that 
$d^{-n}[f^{-n}D^s]\dot\wedge\mu^-$ must be nonzero for some $n$.  This produces
intersections between stable and unstable manifolds.

\proclaim Lemma 9.1.  Let $P$ be a Pesin box, and $\Gamma^u$ be the
corresponding lamination defined in \S 8.   If $p$ is a saddle point, then
$W^s(p)$ must intersect $\lambda^s$ almost every disk of $\Gamma^u$, and the
tangential intersections are an isolated subset of $W^s(p)$.
 
\give Proof.  Let $D\subset W^s(p)$ be a relatively compact open set.  If $D$
contains $p$, then $\mu^-|_D\ne0$, so $d^{-n}[f^{-n}D]\to c\mu^+$ for some
$c>0$ as $n\to\infty$.  Let us suppose that there is a subset $E\subset P^u$
such that the corresponding unstable disks are disjoint from $W^s(p)$.  Let us
define
$$\nu^-_E=\int_{a\in E}\lambda^s(a)\,[\Gamma^u(a)].$$
By Lemma 8.2, it follows that
$\nu^-_E\wedge d^{-n}[f^{-n}D]\to c\nu^-_E\wedge\mu^+$.
By Lemma 8.3, the left
hand side must be zero.  But $\mu^+\ge\nu^+\equiv \mu^+\contract\Gamma^s$, and
$${\bf M}[\nu^-_E\wedge\nu^+]
=\lambda^s(E)\lambda^u(P^s)\ne0.$$
Thus we must have $\lambda^s(E)=0$, which completes the proof of the first
part.  The tangential intersections are isolated by Lemma 6.3.\qed

\proclaim Theorem 9.2.  If $p$ is a saddle point of $f$, then $p\in J^*$.
 
\give Proof.  Let $P$ be a Pesin box.  By Lemma 9.1, there is a compact subset
$K^u\subset W^u(p)$ and a subset $\cG^s_0$ of the leaves of $\Gamma^s$
 such that (i) $\mu^u(\cG^s_0)>0$, (ii) for each $a^u\in K^u$ there
is a leaf $\Gamma^s(x)$ of $\Gamma^s$ such that $\{a^u\}=K^u\cap \Gamma^s(x)$,
and
(iii) the angle of intersection of $\Gamma^s(a^u)$ and $W^u(p)$ is greater than
$\theta_0>0$.  Similarly, we may find subsets $K^s\subset W^s(p)$ and
the family of leaves $\cG^u_0$ with analogous properties.
 
Let us choose a coordinate system in a neighborhood $U$ of $p$ so that $p=0$,
$U=\{|x|<1,\, |y|<1\}$, $\{y=0\}\cap U$ is the component of $W^s(p)\cap U$
containing 0, and $\{x=0\}\cap U$ is the component of $W^u(p)\cap U$ containing
0.  By the Lambda Lemma (see e.g., [PdM]),
 we may take $n$ sufficiently large that
$f^{-n}K^u\subset\{x=0,|y|<\epsilon\}$, and for each $a^u\in K^u$ the portion
of $f^{-n}W^s(a^u)\cap U$ passing through $f^{-n}a^u$ (denoted
$W^s_1(f^{-n}a^u)$), is uniformly $C^1$ close to $\{y=0,|x|<1\}$.   Similarly,
for each $a^s\in K^s$, the portion of $f^nW^u(a^s)\cap U$ passing through
$f^na^s$ (denoted $W^u_1(f^na^s)$), is $C^1$ close to $\{x=0,|y|<1\}$.
 
Thus we may choose $n$ large enough that every pair of
manifolds $W^s_1(f^{-n}a^u)$ and $W^u_1(f^{n}a^s)$ have nonempty intersection.
Finally, since $\lambda^u(K^u)>0$, it follows that
$$\lambda^u_n:=f^{-n}_*(\lambda^u\contract K^u)=\mu^+|_{\{x=0,|y|<1\}\cap
f^{-n}K^u}\ne0.$$
Thus
$$\nu_n^+:=\int_{a\in f^{-n}K^u}\lambda^u_n(a)\,[W^s_1(a)]\ne0.$$
With an analogous definition for $\nu^-_n$, we have $\mu\contract
U\ge\nu_n^+\wedge\nu^-\ne0$.  Thus $p$ is in the support of $\mu$.\qed
 
Combining Theorem 9.2 with the density of saddle points proved in [BS4] gives
the following characterization of $J^*$.
 
\proclaim Corollary 9.3.  $J^*$ is the closure of the set of saddle points.
 
\proclaim Corollary 9.4. Any hyperbolic measure has support contained in $J^*$.
 
\give Proof. According to Katok [K, Theorem 8], periodic saddle points are dense
in the support of any hyperbolic measure. Thus the Corollary follows from
 Theorem 9.2.\qed
 
\proclaim Corollary 9.5. If  $p$ is a saddle point for $f$, then
every transverse intersection of $W^s(p)$ and $W^u(p)$ is in $J^*$.
 
\give Proof. By the Birkhoff-Smale theorem every transverse homoclinic
intersection is the limit of saddle points. So the corollary follows from
Theorem 9.2 and the fact that $J^*$ is closed.\qed
 
\proclaim Theorem 9.6. If $p$ and $q$ are saddle points for $f$, then
the set of transverse intersections of $W^s(p)$ and $W^u(q)$ is dense in $J^*$.
 
\give Proof.  Let $U$ be an open set with $U\cap J^*\ne\emptyset$.  Then there
exists a Pesin box $P\subset U$.  Every pair of points $x_1,x_2\in P$
 have the property that $W^s_r(x_1)$ intersects $W^u_r(x_2)$
in a unique point in $P$ and the intersection is transverse.  Further, there
exists $\epsilon>0$ such that for any smooth manifolds $M^s$ and $M^u$ such that
${\rm dist}_{C^1}(M^s,W^s_r(x))<\epsilon$ and
${\rm dist}_{C^1}(M^u,W^u_r(x))<\epsilon$,
then $M^s$ intersects $M^u$ in a unique point in $U$ and the intersection is
transverse.   By Lemma 9.1, there exist $x_1,x_2\in P$ and $y_1\in W^s(p)\cap
W^u_r(x_1)$ and $y_2\in W^u(q)\cap W^s_r(x_2)$.  For infinitely many values
$n_j\to\infty$ we have $f^{-n_j}x_1\in P$.  By the Lambda Lemma, there is a
portion $D^s_j$ of $f^{-n_j}W^s(p)$ which contains $f^{-n_j}y_1$ and lies as a
graph over $W^s_r(f^{-n_j}x_1)$.  Further, $D^s_j$ approaches
$W^s_r(f^{-n_j}x_1)$ in $C^1$.  Similarly, $f^{m_k}y_2\in P$ for infinitely many
$m_k\to\infty$, and there is a portion $D^u_k\subset f^{m_k}W^u(q)=W^u(q)$ which
approaches $W^u_r(f^{m_j}x_2)$ in $C^1$.  For $j,k$ large, this $C^1$ distance
is less than $\epsilon$, and it follows that $D^s_j$ and $D^u_k$ have a
transverse intersection.  Since $W^s(p)\cap W^u(q)\cap U\supset D^s_j\cap
D^u_k\cap U\ne\emptyset$, it follows that the set of transverse intersections
of $W^s(p)$ and $W^u(q)$ is dense in $J^*$.\qed
 
A saddle point $p$ is said to generate a homoclinic intersection
 if $W^s(p)$ and $W^u(q)$ intersect in points other than $p$.
 
\proclaim Corollary 9.7. Every saddle point generates a  homoclinic
intersection.
 
\proclaim Proposition 9.8.  If  $p$ is a saddle point for $f$, then
every point in the intersection of $W^s(p)$ and $W^u(q)$ is a limit of
transverse intersections of $W^s(p)$ and $W^u(q)$.
 
\give Proof. Let $x\in W^s(p)\cap W^u(q)$. Choose a coordinate system at
$p$ as in the proof of Theorem 9.2. Replacing $x$ by $f^n(x)$ we may assume
 that
$x\in \{y=0\}\cap U$. Now by Theorem 9.4 $W^u(p)$ has a transverse
intersection with $W^s(p)$. The Lambda Lemma implies that there are
components of $W^u(p)\cap U$ arbitrarily close to $\{y=0\}\cap U$. Lemma 6.4
then completes the proof.\qed

\proclaim Theorem 9.9. If  $p$ is a saddle point for $f$, then
every intersection of $W^s(p)$ and $W^u(q)$ is in $J^*$.
 
\give Proof. This follows from Theorem 9.6 and the previous proposition.\qed
\smallbreak

\section 10.  Applications to Real Henon mappings
 
Consider a polynomial diffeomorphism $f$ with real coefficients.
$f$ leaves invariant the real subspace $\real2$. In this section we will denote
$f:\cx2\to\cx2$ by $f_\C$ and $f|_{\real2}$ by $f_\R$. Recall that if $d>1$,
then $d$ can be defined as the  minimal degree of any map conjugate to $f$. This
number can be computed from $f_\R$ without making reference to $f_\C$.
In [FM] it is shown that $h_{top}(f_\R)\le \log d$. In this section we
investigate maps for which equality holds.
 
 Hyperbolic ``$d$-fold'' horseshoes (see [FM] and [HO]) are examples of maps of
maximal entropy but these are not the only examples. The set of horseshoes is
 open
in parameter space and continuity of the entropy function ([Mi]) shows the set of
maps of entropy $\log d$ is closed. Thus the set of parameters of maps with
 entropy
$\log d$ contains the closure of the set of horseshoe parameters. It would be
interesting to know whether it contains other maps as well.
 
\proclaim Theorem 10.1. The following are equivalent:
\item{(1)} $h_{top}(f_{\R})=\log d$
\item{(2)} $\mu(\real2)>0$
\item{(3)} $J^*\subset\real2$
\item{(4)} $K\subset\real2$
\item{(5)} Every periodic point of $f_{\C}$ is contained in $\real2$
\item{(6)} If $p$ is a periodic saddle point then
$W^s(p,f_\C)\cap W^u(p,f_\C)$ is contained in $\real2$.
\vskip0pt
\noindent Any of these conditions implies: 
\item{(7)} $J^*=J=K$.

\give Proof. The result of
Newhouse, (Theorem 2.2), shows that $f_{\R}$ possesses a probability measure
 $\nu$ of
maximal entropy. If (1) holds then the entropy of $\nu$ is $\log d$ so $\nu=\mu$
 by the
uniqueness result, Theorem 3.1. Thus (1) implies  $\mu(\real2)=1$.
 
Assume that (2) holds. Since $\mu$ is ergodic and $\real2$ is and invariant set
 of
positive measure we have $\mu(\real2)=1$. Since $\real2$ is a closed set of full
measure the support of $\mu$ is contained in $\real2$. But the support of $\mu$
 is
$J^*$ so (2) implies (3).
 
 We will show
that if $J^*$ is real then $K$ is real. Recall that $J^*$ is the Shilov boundary
 of
$K$ which is the minimal closed set $S\subset K$ with the property that for any
polynomial $P$ the maximum value of $|P|$ on $K$ is equal to the maximum value
 of
$|P|$ on $S$. It is a general fact that if the Shilov boundary of a set is real
then the set is real. We recall the proof. Assume that $K$ is not real. Say
$p=(z_1,z_2)$ is in $K$ but not in $\real2$. Either $z_1$ or $z_2$ is not real.
For definiteness assume that $z_1\notin\real{}$. Let $J_1$ be the projection of
$J^*$ onto the first coordinate. Runge's theorem assures the existence of a
complex polynomial $P_1(z)$ so that $|P_1(J_1)|<1/10$ and $|P_1(z_1)|>1$. Thus
 the
polynomial $P(z_1,z_2)=P_1(z_1)$ takes its maximum value outside of $J^*$
contradicting our assumption. Thus (3) implies (4).
 
If (4) holds then $\mu$ is supported in $\real2$ so $f_{\R}$ has a measure of
entropy $\log d$ so (4) implies (1). This demonstrates the equivalence of
conditions (1) through (4).
 
We prove the equivalence of (5).  Since every periodic point is in $K$, (4)
implies (5). Since periodic points are dense in $J^*$, ([BS3] Theorem 3.4), we
have (5) implies (3).
 
We prove the equivalence of (6).  Since every point in $W^s(p,f_\C)\cap
 W^u(p,f_\C)$
has a bounded orbit this set is contained in $K$, thus (4) implies (6). Since
$W^s(p,f_\C)\cap W^u(p,f_\C)$ is dense in $J^*$, (Corollary 9.3), we have (6)
implies (3).
 
To show that these conditions imply  (7) we argue as follows.
By (4) $K\subset\real2$. The Stone-Weierstrass theorem implies that any
 continuous
function on $K\subset\real2$ can be approximated by a polynomial function. This
implies that the Shilov boundary of $K$ is all of $K$. Thus $J^*=K$. Since
$J^*\subset J\subset K$ we have $J^*=J=K$ .\qed
 
The following result gives some consequences of the equivalent conditions
described in Theorem 4.1.
Note that all these results can be {\it stated} in terms of $f_{\R}$ without
 reference
to $f_\C$ or $\cx2$. Nevetheless our {\it proofs} of these results
require complex techniques.
 
\proclaim Theorem 10.2. Let $f_{\R}$ be a polynomial diffeomorphism of $\real2$
 with
entropy $\log d$ then:
\item{(1)} $f_\R$ has a unique measure of maximal entropy
\item{(2)} $f_{\R}|K$ is topologically mixing
\item{(3)} $f_{\R}$ has no sinks
\item{(4)} Periodic points are dense in the set of bounded orbits
\item{(5)} For any periodic saddle point $W^s(f_{\R},p)\cap W^u(f_{\R},p)$
is dense in the set of bounded orbits.
 
\give Proof. Since $f_\C$ has a unique measure of maximal entropy it follows
 that
$f_\R$ has a unique measure of maximal entropy when $h(f_\R)=h(f_\C)$.
 
We prove (2). By [BS3] $f$ is mixing for the measure $\mu$. It follows that $f$
is topologically mixing on the support of $\mu$ which is $J^*$. By assertion (7)
 of
Theorem 4.1 $J^*=K$.
 
 To prove (3) we note that a sink orbit
is in $K$ but not in $J$ (a sink is in the interior of $K^+$ but $J=\bd
 K^+\cap\bd
K^-$. So (1) implies that $f$ has no sinks. In the volume
preserving case the same argument shows that $f$ has no
linearizable elliptic points.
 
 Assertion (4) follows from assertion (7) of Theorem 10.1
because periodic points are dense in $J^*$ and the set of bounded orbits is $K$.
 
Assertion (5) follows from (1) because homoclinic intersections
are dense in $J^*$ (Corollary 9.3) and $J^*=K$ from assertion (7) of Theorem
 4.1.
This completes the proof of the theorem. \qed
 
\give Remark.  Let us mention the real quadratic mapping
 $h:\RR^2\mapsto\RR^2$ given
by $(x,y)\to(1-ax^2+y,bx)$ with $a=1.4$ and $b=.3$, which was considered in
detail by H\'enon.  There are eight solutions of $\{(x,y)\in \cx2:h^3_{\bf
C}(x,y)=(x,y)\}$, two of which are  real fixed points, and the other six lie
in two cycles of period 3.  Numerical computation suggests that the 3-cycles
consist of nonreal points. Paul Pedersen gave a mathematical proof that this
is indeed the case ([P]). It follows from Theorem 10.1 that $h$ has entropy
 strictly less than $\log 2$.  And for any saddle $p$, $W^s(p)\cap W^u(p)$
contains points outside $\RR^2$.

\section 11.  Appendix: Concluding Remarks
 
 Some remarks on the logical inter-relationships between the various sections
of this paper are in order.  This paper was organized so that the first
methods used were Pesin theory and entropy; and the first main results
obtained were the identification of the conditional measures and the
uniqueness of the measure of maximal entropy.  The logical progression we
have adopted was not the only one possible.  What follows is an outline of
a different order in which arguments from this paper could be presented.  In
this scenario, the use of entropy comes only at the end.  And this
organization leads to new proofs both that the entropy of $\mu$ is $\log d$
and that the topological entropy of $f$ is $\log d$. 
 \smallskip
 {\sl Step 1.} Let $P$ be a Pesin box, $\Gamma^s$ be the corresponding stable
lamination as defined in \S8. First prove the holonomy invariance of measures
induced by $\mu^+$ along this lamination (Lemma 4.4).
Thus $\mu^+$ induces a transversal measure on $\Gamma^s$. \smallskip
 
{\sl Step 2.} From Lemma 4.4 we deduce that $\mu^+\contract\Gamma^s$ is a
uniformly laminar current (Lemma 8.1).  Next we prove  Lemmas 8.2 and 8.3, and
it follows that $$\mu\contract P=(\mu^+\contract\Gamma^s)\dot
\wedge(\mu^-\contract\Gamma^u).$$  Thus $\mu\contract P$ has a product
structure.\smallskip
 
{\sl Step 3.} The product structure of $\mu$ on sets $P$ of positive
measure, and especially the fact that the wedge product is equal to the
intersection $\dot\wedge$, allows us to prove the results of \S9 concerning
saddle points.\smallskip
 
{\sl Step 4.} The product structure of $\mu$ on the sets $P$ also implies that
the conditional measures of $\mu$ on the unstable leaves are induced by $\mu^+$. 
\smallskip

{\sl Step 5.}  Up to this point, these arguments have not involved entropy.  But
now the facts obtained in the previous steps may be used to calculate the entropy
of $\mu$ via the formula: $$h_\mu(f)=\int\log J_\mu^u\;\mu=\log d.$$\smallskip
 
{\sl Step 6.}  The argument in \S3 shows that $h_\nu(f)<\log d$ for any
invariant measure $\nu\ne\mu$. Hence, $\mu$ is the unique measure of
maximal entropy.
By the Variational Principle, this also gives us an alternative proof
of the Friedland-Milnor-Smillie formula for the topological entropy:
$h(f)=h_\mu(f)=\log d$.\smallskip

\bigbreak
\centerline{\bf References}
\medskip
\def\and{{\&\ }}

\def\endref{\smallskip}
\def\key#1 {\item{[#1]}}
\def\by#1{{#1,}}
\def\paper#1{{\it #1,}}
\def\jour{}
\def\book#1{{\it #1\ }}
\def\publ#1{\ #1}
\def\vol#1 {\ {\bf#1}}
\def\no#1 {{\ no. #1}}
\def\yr#1 {\ (#1)}
\def\pages#1  {\ #1}

\def\paperinfo#1{{\ #1}}
\def\bookinfo#1{\ #1}
\def\inbook#1{#1}
\def\ed#1 {#1}
\def\moreref{;\ }

\key Be
\by{E. Bedford}
\paper{Iteration of polynomial automorphisms of $\bf C^2$}
\inbook{Proceedings of the I.C.M. 1990, Kyoto, Japan}
\yr 1991
\pages 847--858.
\endref

\key BS1
\by{E. Bedford \& J. Smillie}
\paper{Polynomial diffeomorphisms of $\bf C^2$: Currents,
       equilibrium measure and hyperbolicity}
\jour Inventiones math.
\vol 87
\yr 1990
\pages 69--99.
\endref

\key BS2
\by{E. Bedford \& J. Smillie}
\paper{Fatou--Bieberbach domains arising from polynomial
    automorphisms}
\jour Indiana U. Math. J.
\vol 40
\yr 1991
\pages 789--792.
\endref

\key BS3
\by{E. Bedford \& J. Smillie}
\paper{Polynomial diffeomorphisms of $\bf C^2$  II: Stable manifolds
      and recurrence}
\jour Jour. AMS
\vol 4
\no 4
\yr 1991
\pages 657--679.
\endref

\key BS4
\by{E. Bedford \& J. Smillie}
\paper{Polynomial diffeomorphims of $\bf C^2$ III: Ergodicity,
      exponents and entropy of the equilibrium measure}
\jour Math. Annalen
\endref

\key BT1
\by{E. Bedford \& B. A. Taylor}
\paper{Fine topology, \u Shilov boundary and $(dd^c)^n$}
\jour J. of Functional Analysis
\vol 72
\yr 1987
\pages 225--251.
\endref

\key BT2
\by{E. Bedford \& B. A. Taylor}
\paper{Positive, closed currents and complex cycles}
\paperinfo{unpublished manuscript}
\yr 1988
\endref
 
\key Bo
R. Bowen,  {\it Equilibrium states and the ergodic theory of Anosov
diffeomorphisms,} Springer Lecture Notes 470, (1975).\smallskip
 
\key Br
\by{H. Brolin}
\paper{Invariant sets under iteration of rational functions}
\jour Ark. Mat.
\vol 6
\yr 1965
\pages 103--144.
\endref

\key C
\by{L. Carleson}
\paper{Complex dynamics}
\paperinfo{U.C.L.A. course notes}
\yr 1990
\endref
 
\key CFS
I.P. Cornfeld, S.V. Fomin, \& Ya.G. Sinai, {\it Ergodic Theory}, Springer-Verlag,
1982.\smallskip
 
\key D
\by{J.--P. Demailly}
\paper{Courants positifs extr\^emaux et conjecture de Hodge}
\jour Invent. Math.
\vol 69
\yr 1982
\pages 347--374.
\endref
 
\key EL
\by{A. Eremenko \& M. Lyubich}
\paper{Dynamics of analytic transformations}
\jour Leningrad J. Math.
\vol 1
\yr 1989
\pages 1--70.
\endref
 
\key FHY
\by{A. Fathi, M. Herman, \& J.--C. Yoccoz}
\paper{A proof of Pesin's stable manifold theorem}
\inbook  in Geometric Dynamics
\ed J. Palis
\publ Springer Lecture Notes 1007
\yr 1983
\pages 177--215.
\endref

\key FS
\by{J.--E. Forn\ae ss and N. Sibony}
\paper{Complex H\'enon mappings in $\bf C^2$ and Fatou Bieberbach
     domains}
\jour Duke Math.\ J.\ 
\vol 65
\yr 1992
\pages 345--380.
\endref

\key FM
\by{S. Friedland \& J. Milnor}
\paper{Dynamical properties of plane polynomial automorphisms}
\jour Ergodic Theory \& Dyn. Syst.
\vol 9
\yr 1989
\pages 67--99.
\endref
 
\key G
\by{M. Gromov}
\paper{ On the entropy of holomorphic maps}
\paperinfo{Unpublished manuscript}
\yr 1980
\endref

\key Ha
W. Hayman, {\it Meromorphic Functions}, Oxford, 1964. \smallskip

\key H
J.H. Hubbard, {\it The H\'enon mapping in the complex domain}, In: M. Barnsley
and S. Demko (eds.), Chaotic Dynamics and Fractals. Academic Press, New York
1986, pp. 101--111.\smallskip

\key HO
\by{J. H. Hubbard \& R. Oberste--Vorth}
\paper{H\'enon mappings in the complex domain}
\paperinfo{in preparation}.
\endref

\key HP
\by{J. H. Hubbard \& P. Papadopol}
\paper{Superattractive fixed points in $\bf C^n$}
\paperinfo{preprint}.
\endref

\key K
\by{A. Katok}
\paper{Lyapunov exponents, entropy and periodic orbits for
           diffeomorphisms}
\jour Pub. Math. de I.H.E.S.
\vol 51
\yr 1980
\pages 137--174.
\endref

\key Kl
\by{M. Klimek}
\book{Pluripotential theory}
\publ Oxford
\yr 1991
\endref
 
\key Kr
W. Krieger, {\it On entropy and generators of measure-preserving
transformations}, Trans. AMS. {\bf 119} (1970), pp 453--464.\smallskip

\key Le
\by{F. Ledrappier}
\paper{Some properties of absolutely continuous invariant measures
       on an interval}
\jour Ergod. Th. \& Dynam. Sys.
\vol 1
\yr 1981
\pages 77-93
\endref

\key LS
\by{F. Ledrappier \& J.--M. Strelcyn}
\paper{A proof of the estimation from below in Pesin's entropy
      formula}
\jour Ergod. Th. \& Dynam. Sys.
\vol 2
\yr 1982
\pages 203--219.
\endref

\key  LY
\by{F. Ledrappier \& L.--S. Young}
\paper{The metric entropy of diffeomorphisms, I}
\jour Annals of Math.
\vol 122
\yr 1985
\pages 509--539
\moreref
\paper{II}
\jour Annals of Math.
\vol 122
\yr 1985
\pages 540--574.
\endref

\key L
\by{P. Lelong}
\book{El\'ements extr\^emaux sur le c\^one des courants positifs
      ferm\'es}
\bookinfo{S\'eminaire P. Lelong (Analyse), Ann\'ee, 1971--1972,
     Lecture Notes in Math.}
\publ Springer
\vol 332
\yr 1972
\endref
 
\key Lyu1
M. Lyubich, {\it Entropy of analytic endomorphisms of the Riemannian sphere},
Funct. Anal. Appl. {\bf 15}(1981), 300--302.\smallskip

\key Lyu2
\by{M. Lyubich}
\paper{Entropy properties of rational endomorphisms of the Riemann
      sphere}
\jour Ergod. Th. \& Dynam. Sys.
\vol 3
\yr 1983
\pages 351--385.
\endref

\key Ma
\by{R Ma\~{n}\'{e}}
\paper{On the uniqueness of the maximizing measure for rational maps}
\jour Bol. Soc. Brasil Math.
\vol 14
\yr 1983
\pages 27--43.
\endref 

\key Mak
\by{N. G. Makarov}
\paper{On the distortion of boundary sets under conformal mappings}
\jour Proc. London Math. Soc.
\vol 51
\yr 1985
\pages 369-384
\endref

 
\key Man
A. Manning, {\it The dimension of the maximal measure for a polynomial map},
Ann. Math. {\bf 119} (1984), 425--430.\smallskip
 
\key Mi
J. Milnor, {\it Non-expansive H\'enon maps} Adv. Math.. 69 (1988),
109--114.\smallskip

\key Mo
\by{A. P. Morse}
\paper{Perfect blankets}
\jour Trans. AMS
\yr 1947
\pages 418--442.
\endref
 
\key N
R. Nevanlinna, {\it Analytic Functions}, Springer-Verlag.\smallskip

\key Ne
S. Newhouse, {\it Continuity properties of entropy}, Annals of Math., {\bf 129}
(1989), 215--235.\smallskip
 
\key OW
\by{D. Ornstein and B. Weiss}
\paper{Statistical properties of chaotic systems}
\jour Bull. AMS
\vol 24
\yr 1991
\pages 11--116.
\endref
 
\key PdM
J. Palis and W. de Melo, {\it Geometric Theory of Dynamical Systems},
Springer-Verlag, New York, 1982.\smallskip
 
\key Pe 
P. Pedersen, personal communication.\smallskip

\key P1
\by{Y. Pesin}
\paper{Characteristic Lyapunov exponents and smooth ergodic theory}
\jour Russian Math. Surveys
\vol 32
\yr 1977
\pages 55--114.
\endref
 
\key P2
\by {Y. Pesin}
\paper{Description of the $\pi$-partition of a diffeomorphism with
           invariant smooth measure}
\jour Mat. Zametki
\vol 21
\yr 1977
\pages 29--44.
\endref

\key Pr
\by{F. Przytycki}
\paper{Hausdorff dimension of harmonic measure on the boundary of an attractive
basin for a holomorphic map} \jour Invent. math.
\vol 80
\yr 1985
\pages 161--179.
\endref

\key PS
\by{C. Pugh \& M. Shub}
\paper{Ergodic attractors}
\jour Trans. AMS
\vol 312
\yr 1989
\pages 1--54.
\endref

\key Ro1
\by{V. A. Rokhlin}
\paper{On the fundamental ideas of measure theory}
\jour Mat. Sbornik 25 (67) 107-150 (1949)
\publ{AMS}
\yr 1952
\endref
 
\key Ro2
V.A. Rohlin, {\it Lectures on entropy theory of measure-preserving
transformation} Russian Math. Surveys, {\bf 22} (1967),  \smallskip
 
\key R1
D. Ruelle, {\it An inequality for the entropy of differential maps,} Bol. Soc.
Brasil Mat. {\bf 9} (1978), 83--87.\smallskip

\key R2
\by{D. Ruelle}
\paper{Ergodic theory of differentiable dynamical systems}
\jour Publ. Math. de I.H.E.S.
\vol 50
\yr 1979
\pages 275--306.
\endref

\key RS
\by{D. Ruelle and D. Sullivan}
\paper{Currents, flows, and diffeomorphisms}
\jour Topology
\vol 14
\yr 1975
\pages 319--327.
\endref
 
\key Si
N. Sibony, Course at UCLA, unpublished manuscript. \smallskip

\key S
\by{J. Smillie}
\paper{The entropy of polynomial diffeomorphisms of $\bf C^2$}
\jour Ergodic Theory \& Dynamical Systems
\vol 10
\yr 1990
\pages 823--827.
\endref

\key Su
\by{D. Sullivan}
\paper{Dynamical study of foliated and complex manifolds}
\jour Invent. Math.
\vol 36
\yr 1976
\pages 225--255.
\endref

\key T
\by{P. Tortrat}
\paper{Aspects potentialistes de l'it\'eration des polyn\^omes}
\jour S\'eminaire de Th\'eorie du Potentiel Paris, No. 8,
\yr 1987,
\publ Springer Lecture Notes 1235.
\endref
 
 
\key W
\by{H. Wu}
\paper{Complex stable manifolds of holomorphic diffeomorphisms
       and complex H\'{e}non mapping}
\paperinfo{ Preprint}
\endref

\key Yg
\by{L.--S. Young}
\paper{Dimension, entropy and Lyapunov exponents}
\jour Ergod. Th. \& Dynam. Sys.
\vol 2
\yr 1982
\pages 109--124.
\endref

\bigskip
\noindent Indiana University, Bloomington, IN 47405
\medskip
\noindent SUNY at Stony Brook, Stony Brook, NY 11794
\medskip
\noindent Cornell University, Ithaca, NY 14853

\bye